\documentclass[12pt, times]{article} 
\usepackage[letterpaper, margin=1in]{geometry}

\usepackage{enumitem}

\usepackage[gen]{eurosym}
\usepackage{amssymb,graphics,amsmath,amsthm,amsopn,amstext,amsfonts}
\newtheorem{theorem}{Theorem}
\newtheorem{proposition}[theorem]{Proposition}
\newtheorem{lemma}[theorem]{Lemma}
\newtheorem{corollary}[theorem]{Corollary}
\newtheorem{definition}[theorem]{Definition}
\newtheorem{example}[theorem]{Example}
\newtheorem{remark}[theorem]{Remark}
\usepackage{epstopdf}
\usepackage{framed}
\usepackage{color}
\usepackage[dvips]{epsfig}
\usepackage{graphicx}
\usepackage{subcaption,multirow,array}
\usepackage{amsmath}
\usepackage{mathtools}
\usepackage{xcolor}
\usepackage{hyperref}
\hypersetup{
	colorlinks= true,
	citecolor=blue,
	linkcolor=blue  
}

\usepackage{color}

\usepackage[utf8]{inputenc}
\usepackage[english]{babel}
\usepackage{endnotes}
\usepackage{changepage}
\usepackage{comment}

\usepackage{cancel}

\newcommand{\argmin}{\operatornamewithlimits{\arg\min}}

\newcolumntype{A}{>{\centering\arraybackslash}p{0.1\textwidth}}
\newcolumntype{B}{>{\centering\arraybackslash}p{0.16\textwidth}}
\newcolumntype{C}{>{\centering\arraybackslash}p{0.08\textwidth}}

\newcommand{\kpd}{k, \boldsymbol p^0, \boldsymbol \delta}
\newcommand{\pdi}{\mathcal P^i_{\delta_i}}

\let\footnote=\endnote

\usepackage{natbib}
 \bibpunct[, ]{(}{)}{,}{a}{}{,}%
\begin{document}

\title{Price Optimization with Practical Constraints}


\author{Xiaojie Wang, Hsin-Chan Huang, Lanshan Han and Alvin Lim\footnote{Corresponding author, alvin.lim@emory.edu}\\
Nielsen Precima, LLC}



\maketitle

%

\begin{abstract}%
In this paper, we study a retailer price optimization problem which includes the practical constraints: maximum number of price changes and minimum amount of price change (if a change is recommended). We provide a closed-form formula for the Euclidean projection onto the feasible set defined by these two constraints, based on which a simple gradient projection algorithm is proposed to solve the price optimization problem. We study the convergence and solution quality of the proposed algorithm. We extend the base model to include upper/lower bounds on the individual product prices and solve it with some adjustments to the gradient projection algorithm. Numerical results are reported to demonstrate the performance of the proposed algorithm.
\end{abstract}%

\section{Introduction}
Price optimization is a classical problem in the business operation of various product and service industries. It is well-known from economics that the demand of a product is affected by its price. Typically, for a consumer product, a lower price for the product leads to more demand for it. However, when multiple products are offered, the situation is often more complex due to the cross-product relationships. The products can be substitutes or complements of each other, and hence the demand on a product depends not only on its own price but also the prices of related products. We often use demand functions to quantify how demands of products are affected by their prices.

When price optimization is applied in various industries under different settings, different sets of business rules need to be imposed to guarantee that optimal prices are reasonable and practical to implement. In this paper, we are particularly interested in finding the optimal prices for a set of products in a brick-and-mortar grocer to maximize profits. In grocery retailing, certain specific business considerations need to be taken into account in setting prices of products. First, prices change frequently to reflect commodity cost changes, supplier incentives, response to consumer behavior, response to market conditions, etc., and price optimization is often implemented to obtain product shelf prices for various store locations. While a growing adoption of electronic shelf tags has been observed recently \citep{bhutani2018electronic}, the use of paper tags that are implemented using limited manual labor is still the norm.  Due to the limited (and often expensive) resources, grocers only update in each week the prices of a limited number of products -- those that have the most significant impact on revenues and profits. Second, while demand functions are typically continuous function defined in the nonnegative orthant and hence can describe the demand change expected from an arbitrary small change in price, consumers may not be sensitive enough to react to extremely small changes in price.  So, when a price decrease is too small, the expected increase in demand for a product may not materialize.  On the other hand, if the price increase is too small, while the demand will likely stay constant (leading to better profits), the cost of implementing the price increase may not justify the return. Therefore, a minimum price change is often imposed on products that are candidates for price changes. These two practical considerations significantly complicate the price optimization problem. In particular, they implicitly introduce certain combinatorial structures into the traditional price optimization problem, which is typically a continuous optimization problem (often also convex under reasonable assumptions) \citep{gallego2019revenue}. A mixed integer programming formulation is, therefore, needed to capture these two practical considerations. This type of optimization problem is notoriously challenging to solve to global optimality \citep{wolsey2007mixed}.

This paper represents a first effort to include the aforementioned practical considerations into price optimization models. We explore various mathematical properties of the resulting formulation and investigate the use of a gradient projection algorithm to find a solution of the optimization problem. The algorithm is demonstrated to produce near-optimal solutions both theoretically and empirically. Moreover, we hope that the proposed approach may be applied to similar problems, such as marketing mix optimization \citep{huang2021marketing} and portfolio optimization \citep{bertsimas2009algorithm}, when analogous practical constraints are imposed.

The rest of the paper is organized as follows. In Section~\ref{sec:review}, we review existing research on price optimization. In Section~\ref{sec:problem}, we provide mathematical models of the price optimization problem under investigation and present preliminaries regarding the mathematical properties of the models. In Section~\ref{sec:algorithm}, we present a gradient projection algorithm to solve the price optimization problems as well as analyses on various properties of the algorithm such as convergence and performance bounds. Numerical results on randomly generated data sets as well as real-world retailer data sets are presented in Section~\ref{sec:numerical_results}. We conclude the paper with a summary of our contributions and a discussion on possible future research directions in Section~\ref{sec:conclusion}.

\section{Literature Review }\label{sec:review}
In this section, we review the state of the art in retail price optimization under various settings. According to how the pricing decisions are made along the temporal dimension, the literature can be divided into ``dynamic'' pricing and ``static'' pricing categories. Dynamic pricing refers to the setting of time-varying prices over a selling period, and prices are typically set based on the changing expected availability of and expected demand on the product over time. Static pricing, on the other hand, refers to the setting of fixed prices over a selling period of products that are continuously replenished and virtually non-perishable \citep{kunz2014demand}. While dynamic pricing has seen applications in retail via e-commerce \citep{boyd2003revenue}, we will focus on static pricing, which is the predominant practice among brick-and-mortar retailers.

Although more research efforts have been invested into dynamic pricing \citep{bitran2003overview, chiang2007overview, goensch2009dynamic, elmaghraby2003dynamic, weatherford1992taxonomy, mcgill1999revenue, gallego2004managing, liu2008choice, bront2009column}, static pricing dominates retail pricing applications \citep{kunz2014demand}, as frequent price changes are not practical to implement at most brick-and-mortar retailers. Even at online retailers, where one may consider frequent price changes come at negligible costs, frequent changes in price may lead consumers' negative perception of an online retailer's brand image, and therefore not desirable for most retailers, particularly for fashion retailers \citep{ferreira2016analytics,levy2004emerging}. However, despite the prevalence of static pricing, limited studies on the static pricing problem exists in the operations research and econometric literature \citep{kunz2014demand}.

\cite{kunz2014demand} observed the research scarcity and importance of static retail pricing, including the case of grocery pricing. 
\cite{levy2004emerging} examined the practice of retail pricing and found seven critical elements that are essential for success in retail pricing: price sensitivity effects, substitution effects, dynamic effects of price promotions, segment-based pricing effects, cross-category effects, retailer costs and discounts, and retail competition.  Several additional factors should also be considered in the implementation of retail pricing: market factors, grouping of products (category management), continuous learning, psychological price thresholds and reference effects, costs of changing prices, and the quality of data.

The classical method in static retail pricing employs rule-based approaches in category management, which is still the prevailing method used in the industry because of the lack of quality data and models to describe the complex interactions among products as well as the time and expertise to conduct in-depth analyses \citep{hall2003category}.
Some studies on static retail pricing in the literature are concerned with ``multi-product pricing'', which focus on understanding interactions and substitutions among products, resulting in the formulation of the pricing problem as a mixed-integer program, which allows for a wide range of practical constraints but is NP-complete and its efficient solution remains a challenge \citep{dobson1988positioning,oren1987multi}.

A demand model is considered the core of a static retail price optimization system that serves as a link between
the input data and the optimization model \citep{garrett2005}. Types of demand models \citep{kunz2014demand} include:
\begin{itemize}[leftmargin=.5in]
    \item absolute: linear, multiplicative, power, exponential, double-log, etc.;
    \item relative: multiple competitive interaction (MCI) and discrete choice model; and
    \item other forms: reference price, game theory, and hedonic pricing.
\end{itemize}
Many efforts were invested in discrete choice models to account for customer choice behaviors. A comprehensive review of choice-based pricing and revenue management can be found in \cite{strauss2018review}.  Discrete choice models are implemented using multinomial logit (MNL) \citep{song2007demand,dong2009dynamic}, general attraction \citep{keller2014efficient}, nested logit \citep{huh2015pricing}, exponomial choice \citep{alptekinouglu2016exponomial}, generalized extreme value (GEV) \citep{zhang2018multiproduct}, and non-parameteric \citep{farias2013nonparametric,rusmevichientong2006nonparametric} and other techniques.
A more comprehensive review of demand functions in decision modeling can be found in \cite{huang2013demand}. We note that static retail pricing is mainly motivated by demand interdependencies and the consideration of cross price effects is therefore fundamental \citep{kunz2014demand}.

Estimation of the unknown parameters in the demand model is essential, because models without quality input data and parameter estimates have little value. In grocery retailing, the challenge of estimating a large number of parameters in a demand model likely contributed to the limited implementation of data-driven static pricing until the recent prominence in the adoption of data analytics. Bayesian estimation models are used to predict sales as a function of price and other factors \citep{kalyanam1996pricing,seetharaman2005models}. An application at Rue La La, an online fashion retailer, uses machine learning techniques (bagged regression tree) to understand historical sales and predict future demands on changing prices \citep{ferreira2016analytics}. With the predicted demand, \cite{ferreira2016analytics} developed an integer programming model to maximize the total revenue, which is difficult to solve in large scale instances, and thus a linear programming bound algorithm was proposed to help achieve optimality.

The parameters of the demand function can also be estimated through online learning. 
A common framework in this area is to obtain information on demands through price experimentation based on some parametric demand choice models, and through experimentation, pricing polices are designed to minimize regret \citep{broder2012dynamic,ferreira2018online, ozer2012oxford}. More examples of demand function modeling in pricing and revenue management applications can be found in \cite{mivsic2020data}.

\section{Problem Statement}\label{sec:problem}
Throughout this paper, we use lower and upper case letters to represent scalars, bold lower case letters to represent vectors, and bold upper case letters to represent matrices. We use superscript to index vectors or matrices and use subscripts to index scalars. All vectors are assumed to be column vectors. We denote the $n$-dimensional real vector space as $\mathbb R^n$, its nonnegative orthant as $\mathbb R_+^n$, and the interior of the nonnegative orthant as $\mathbb R_{++}^n$. For a vector $\boldsymbol v \in \mathbb R^n$, its $i$-th element is denoted by $v_i$ and its transpose is written as $\boldsymbol v^T$. For a matrix $\boldsymbol A \in \mathbb R^{n \times m}$, its $(i,j)$ element is denoted by $a_{ij}$ and its transpose is written as $\boldsymbol A^T$.  We denote the all zero vector and all one vector of length $n$ by $\boldsymbol 0_n$ and $\boldsymbol 1_n$, respectively. We denote the $n\times n$ identity matrix by $\boldsymbol I_n$. The two-norm of a vector $\boldsymbol v \in \mathbb R^n$ is denoted by $\|\boldsymbol v\|_2$, i.e., $\|\boldsymbol v\|_2 = \sqrt{\boldsymbol v^T \boldsymbol v}$. The infinite-norm of a vector $\boldsymbol v \in \mathbb R^n $ is denoted by $\|\boldsymbol v\|_\infty$, i.e., $\|\boldsymbol v\|_\infty = \max_{i=1,\cdots,n}\{|v_i|\}$. We also define the zero-norm of a vector $\boldsymbol v \in \mathbb R^n$, denoted by $\|\boldsymbol v\|_0$, as the cardinality of the set of nonzero elements of $\boldsymbol v$. Given a function $g:\mathbb R^n \mapsto \mathbb R$, we let $\nabla g(\boldsymbol x)$ be its gradient evaluated at $\boldsymbol x$, and let $\nabla_i g(\boldsymbol x)$ be its $i$-th element. Given a matrix $\boldsymbol A \in \mathbb{R}^{n\times n}$ and a subset $\alpha \subseteq \{1,\cdots,n\}$, we let $\boldsymbol A_{\alpha, \bullet}$ be the $|\alpha| \times n$ submatrix formed by the rows of $\boldsymbol A$ with indices in $\alpha$, with $|\alpha|$ being the cardinality of $\alpha$. Similarly we let $\boldsymbol A_{\bullet,\alpha}$ be the $n\times |\alpha|$ submatrix formed by the columns of $\boldsymbol A$ with indices in $\alpha$, and $\boldsymbol A_{\alpha,\alpha}$ be the principal submatrix with column and row indices in $\alpha$. In the special case where $\alpha$ is a singleton, i.e., $\alpha = \{i\}$, we simplify the notation by using $\boldsymbol A_{\bullet, i}$ and $\boldsymbol A_{i,\bullet}$ to represent the $i$th column and $i$th row of $\boldsymbol A$, respectively. For a vector $\boldsymbol v\in \mathbb R^n$, $\boldsymbol v_\alpha$ is the subvector of length $|\alpha|$ composed of elements of $\boldsymbol v$ with indices in $\alpha$.

We consider a set of $n$ products, indexed by $i=1,\cdots,n$. Let $\boldsymbol p= \left(p_i\right)_{i=1}^n$ be the vector of the unit prices of the products, $\boldsymbol c=\left(c_i\right)_{i=1}^n$ be the vector of the unit costs of the products, and $v(\boldsymbol p):\mathbb R^n \mapsto \mathbb R^n$ be the \emph{(volume of) demand function} which estimates the expected volume of demand as a function of the prices of the products. The price optimization problem is to maximize the profit, possibly under various business rules. Mathematically, this problem is written as:
\begin{equation}\label{eq:price_optimization}
\max_{\boldsymbol p \in \mathcal P} \, \, Z(\boldsymbol p) \, \triangleq \, (\boldsymbol p-\boldsymbol c)^T v(\boldsymbol p),
\end{equation}
where $\mathcal P$ represents the set of all price vectors satisfying defined business rules. While there are many possible business rules such as the maximum percentage change from a baseline price, maximum price gap to a competitor's prices, or price relationship requirements for different package sizes of the same product, these rules can typically be written as linear inequality constraints. These rules create a convex set of feasible price vectors (a polyhedron), which is often considered benign in optimization approaches. On the other hand, there are other business rules that cause significant challenges by creating non-convex feasible set of price vector $\mathcal P$. In this paper, we consider two such constraints that are of significant practical relevance:
\begin{itemize}[leftmargin=.2in]
\item \textbf{Maximum number of price changes}. For most  retailers, the price optimization is conducted routinely with the most recent data collected. While it is possible for an online retailer to change the price of many or even all of their products in a very short time period, it is practically impossible for a traditional brick-and-mortar store to do so, as they only have limited resources to update the paper price tags on the shelves. Therefore, it is of crucial importance that we only change the prices of a priority subset of products that can generate the most incremental profits. Let the baseline price vector be $\boldsymbol p^0 = \left(p^0_i\right)_{i=1}^n$. Assuming that the resources available allow changing the prices of at most $k$ products, this requirement can be translated to a constraint of the following format:
    $$\|\boldsymbol p-\boldsymbol p^0\|_0 \, \leq \, k.$$

\item \textbf{Minimum amount of price change}. As mentioned earlier, extremely small price changes may not deliver the expected return in profits due to the consumers insensitivity to small price changes in the case of a price decrease and the cost of implementing the change in the case of a price increase. In addition, many retailers also apply price ending rules as part of their price image strategy, such as prices ending with 9 cents. To address these two issues, in practice we often need to impose a minimum price change requirement. More specifically, product $i$'s price can either stay the same or increase/decrease by a minimum amount $\delta_i \geq 0$. Mathematically, we have:
    $$p_i \, \in \, \{p^0_i\} \cup [p^0_i + \delta_i, \infty) \cup (-\infty, p^0_i - \delta_i] \, \subseteq \, \mathbb R.$$
    We write $\pdi \triangleq \{p^0_i\} \cup [p^0_i + \delta_i, \infty) \cup (-\infty, p^0_i - \delta_i]$ and $\mathcal P_{\boldsymbol \delta} \triangleq \bigotimes_{i=1}^n \pdi \subseteq \mathbb R^n$ with $\boldsymbol \delta \triangleq \left(\delta_i\right)_{i=1}^n$. Sometimes, we also have upper and lower bounds, denoted by $u_i$ and $l_i$, on product $i$'s price, in which case we write $\mathcal P_{\delta_i,l_i,u_i} \triangleq \{p^0_i\} \cup [p^0_i + \delta_i, u_i] \cup [l_i, p^0_i - \delta_i]$ and $\mathcal P_{\boldsymbol \delta, \boldsymbol l, \boldsymbol u} \triangleq \bigotimes_{i=1}^n \mathcal P_{\delta_i,l_i,u_i}$, with the $\boldsymbol u$ and $\boldsymbol l$ being vectors of upper bounds and lower bounds, i.e., $\boldsymbol u \triangleq (u_i)_{i=1}^n$ and $\boldsymbol l \triangleq (l_i)_{i=1}^n$. An illustration of set $\mathcal P_{\delta_i,l_i,u_i}$ is given in Figure \ref{fig:price_set}.
\begin{figure}[h]
\centering
\includegraphics[scale=.25]{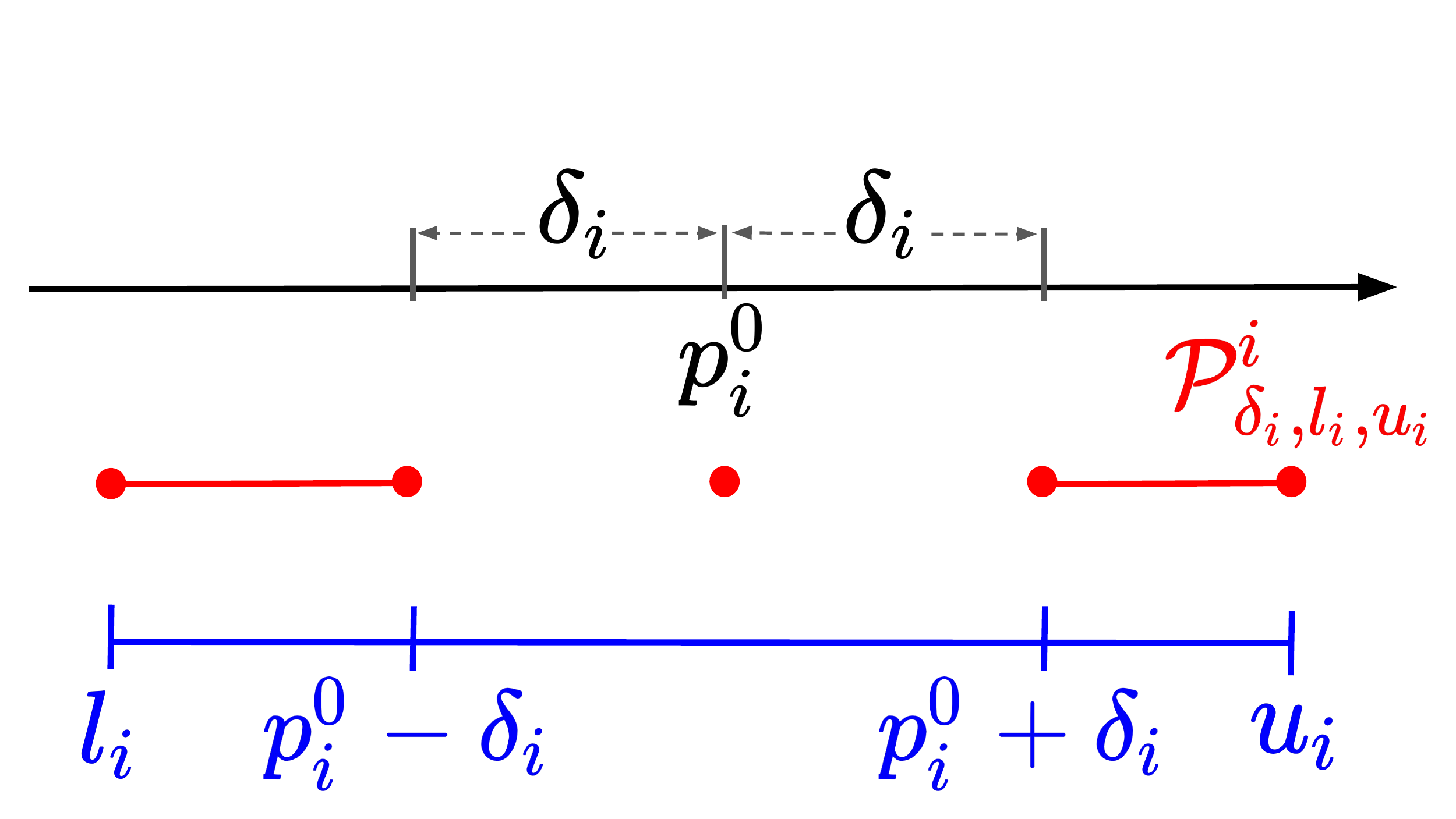}
\caption{ \centering Set of allowed prices for a single product shown in red.}
\label{fig:price_set}
\end{figure}
\end{itemize}
With the above two practical requirements, given baseline price vector $\boldsymbol p^0$ and price change threshold vector $\boldsymbol \delta$, we define $$\Omega_{\kpd} \, \triangleq \, \left\{\boldsymbol p \, \left| \, \|\boldsymbol p- \boldsymbol p^0\|_0 \leq k, \, \boldsymbol p \in \mathcal P_{\boldsymbol \delta}\right. \right\},$$
which is the set of feasible price vectors.  When the upper/lower bounds on individual products exist, this set of feasible price vectors is
$$ \Omega_{k, \boldsymbol p^0, \boldsymbol \delta, \boldsymbol l, \boldsymbol u}\, \triangleq \, \left\{\boldsymbol p \, \left| \, \|\boldsymbol p- \boldsymbol p^0\|_0 \leq k, \, \boldsymbol p \in \mathcal P_{\boldsymbol \delta, \boldsymbol l, \boldsymbol u}\right. \right\}.$$
In this paper, we mainly focus on the case without bounds as follows:
\begin{equation} \label{eq:price_opt_const_base}
\begin{array}{rll}
\displaystyle{ \max_{\boldsymbol p} }&  (\boldsymbol p - \boldsymbol c)^T v(\boldsymbol p) \\[5pt]
\mbox {s.t.}  & \boldsymbol p \, \in \, \Omega_{\kpd},
\end{array}
\end{equation}
but we will also consider the case with upper and lower bounds as a direct extension.

\subsection{Linear Demand Functions}
In this paper, we consider a linear demand function of the form $v(\boldsymbol p) = \boldsymbol a-\boldsymbol D\boldsymbol p$ with $\boldsymbol a \in \mathbb R^n$ and $\boldsymbol D = \left(d_{ij}\right)^n_{i,j=1}\in \mathbb R^{n \times n}$, where for $i =1,\cdots,n$, $-d_{ii}$ represents the effect of product $i$'s own price on its demand $v_i$, and $-d_{ij}$ represents the effect of the cross price $p_j$ of product $j$ on the demand $v_i$ of product $i$. Typically, $\boldsymbol D$ is not symmetric, and the symmetric part of $\boldsymbol D$ is denoted by $\boldsymbol S \triangleq (\boldsymbol D + \boldsymbol D^T)$. Often we make the following two assumptions on $\boldsymbol D$ or $\boldsymbol S$ \citep{gallego2019revenue}.
\begin{itemize}[leftmargin=.5in]
\item[(A1).] Own price effects are all negative and cross price effects are all positive. Mathematically, we have $\boldsymbol S = \left(s_{ij}\right)^n_{i,j=1}$ is positive definite and $s_{ij} \leq 0$ for all $i\neq j$.
\item[(A2).] The demands when prices are all $0$'s are positive, the costs are all positive, and the demands when the prices are at the costs are nonnegative. Mathematically, we have $\boldsymbol a \in \mathbb R_{++}^n$, $\boldsymbol c \in \mathbb R^n_{++}$, and $\boldsymbol a - \boldsymbol D \boldsymbol c \in \mathbb R_+^n$.
\end{itemize}
Note that the cross price effects limitation of (A1) restricts us to consider only the price effects of substitute products.  Under a linear demand function, it is easy to see that $$ Z(\boldsymbol p) \, =  \, (\boldsymbol p - \boldsymbol c)^T(\boldsymbol a-\boldsymbol D\boldsymbol p) \, = \, -\boldsymbol p^T \boldsymbol D \boldsymbol p + \boldsymbol a^T \boldsymbol p + \boldsymbol c^T\boldsymbol D \boldsymbol p - \boldsymbol c^T \boldsymbol a \, = \, -\frac{1}{2}\boldsymbol p^T \boldsymbol S \boldsymbol p + (\boldsymbol a+\boldsymbol D^T \boldsymbol c)^T \boldsymbol p - \boldsymbol c^T \boldsymbol a.$$
Therefore, the base price optimization problem (\ref{eq:price_opt_const_base}) can be written as the following equivalent minimization problem.
\begin{equation} \label{eq:price_opt_const_base_linear_demand}
\begin{array}{rll}
\displaystyle{ \min_{\boldsymbol p} }&  \mathcal Q(\boldsymbol p) \, \triangleq \, \displaystyle{ \frac{1}{2}\boldsymbol p^T \boldsymbol S \boldsymbol p - (\boldsymbol a + \boldsymbol D^T \boldsymbol c)^T \boldsymbol p } \\[5pt]
\mbox {s.t.}  & \boldsymbol p \, \in \, \Omega_{\kpd}.
\end{array}
\end{equation}
Define
\begin{equation}\label{eq:unconst_optimizer}
\widehat{\boldsymbol p} \, \triangleq \, \argmin_{\boldsymbol p} \,\, \mathcal Q(\boldsymbol p),
\end{equation}
i.e., the unconstrained unique minimizer (due to the positive definiteness of $\boldsymbol S$) of the function $\mathcal Q(\boldsymbol p)$. It is easy to see that
$$\widehat{\boldsymbol p} \, = \, \boldsymbol S^{-1}(\boldsymbol a+\boldsymbol D^T \boldsymbol c).$$
Let $\widehat{\mathcal Q} \, \triangleq \,\mathcal Q(\widehat{\boldsymbol p})$, i.e., the unconstrained optimal value. Let $\mathcal Q^*$ be the optimal value of (\ref{eq:price_opt_const_base_linear_demand}). It is clear $\widehat{\mathcal Q} \leq \mathcal Q^*$. Note that the existence and uniqueness of $\widehat{\boldsymbol p}$ is guaranteed by the positive definiteness of $\boldsymbol S$. We have the following lemma, which is known in the literature. We include a short proof for completeness.
\begin{lemma}\label{lm:profitbility}
Assume (A1) and (A2) hold. Let $\widehat{\boldsymbol p}$ be defined by (\ref{eq:unconst_optimizer}), then $$\widehat{\boldsymbol p} \, \geq \, \boldsymbol c.$$
\end{lemma}
\proof{}
It is easy to see that
\begin{eqnarray}
\widehat{\boldsymbol p} - \boldsymbol c & = & \boldsymbol S^{-1}(\boldsymbol a + \boldsymbol D^T \boldsymbol c) - \boldsymbol c
\, = \,  \boldsymbol S^{-1}(\boldsymbol a + \boldsymbol D^T \boldsymbol c - \boldsymbol S \boldsymbol c)\nonumber\\
& = & \boldsymbol S^{-1}(\boldsymbol a + \boldsymbol D^T \boldsymbol c - (\boldsymbol D+\boldsymbol D^T) \boldsymbol c) \, = \, \boldsymbol S^{-1}(\boldsymbol a  - \boldsymbol D \boldsymbol c).\nonumber
\end{eqnarray}
By Assumption (A1), $\boldsymbol S$ is a \emph{Stieltjes} matrix (named after Thomas Joannes Stieltjes), which is necessarily an M-matrix \citep{berman1994nonnegative}. Therefore, $\boldsymbol S$ is always invertible and $\boldsymbol S^{-1}$ is a symmetric nonnegative matrix, i.e., all the elements of $\boldsymbol S^{-1}$ are non-negative. By (A2), $ (\boldsymbol a-\boldsymbol D \boldsymbol c) \geq 0$. Therefore $\widehat{\boldsymbol p} - \boldsymbol c \geq 0$, and we have the desired result.
\endproof
Lemma \ref{lm:profitbility} shows that when we don't have any constraints on the prices, assumptions (A1) and (A2) ensure that the optimal prices are profitable for each product. We call a price vector $\boldsymbol p$ \emph{profitable} if $\boldsymbol p \geq \boldsymbol c$, where the $\geq$ relationship is element-wise. We next study the case of optimization problem (\ref{eq:price_opt_const_base_linear_demand}). We make the following additional assumption:
\begin{itemize}[leftmargin=.5in]
  \item[(A3).] $\boldsymbol p^0-\boldsymbol \delta$ is profitable, i.e., $\boldsymbol p^0 - \boldsymbol \delta \geq \boldsymbol c$.
\end{itemize}
Before we present our main result about profitability, we first prove a technical lemma.

\begin{lemma}\label{lm:partial optimality}
Let $\boldsymbol p^{*}$ be a global optimal solution to (\ref{eq:price_opt_const_base_linear_demand}). Define an index set $\omega^* \, = \, \{ i \, | \, \boldsymbol p^{*}_i = p^0_i\} \cup \{ i \, | \, \boldsymbol p^{*}_i = p^0_i+\delta_i\} \cup \{ i \, | \, \boldsymbol p^{*}_i = p^0_i-\delta_i\},$ and its complement set $\overline {\omega^*} = \{ i \, | \, i \notin \omega^*\}$. Then
\begin{equation}\label{eq:partial_optimial}
\boldsymbol p^{*}_{\overline {\omega^*}} \, = \, \boldsymbol S^{-1}_{\overline {\omega^*},\overline {\omega^*}}\left( \boldsymbol a_{\overline {\omega^*}}+ \boldsymbol D_{\bullet, \overline {\omega^*}}^T \boldsymbol c -\boldsymbol S_{\overline {\omega^*},\omega^*} \boldsymbol p^*_{\omega^*}\right).
\end{equation}
\end{lemma}
\proof{}
We claim that $\boldsymbol p^{*}_{\overline {\omega^*}}$ necessarily minimize the following quadratic function of $\boldsymbol p_{\overline{\omega^*}} $:
\begin{equation}\label{eq:restricted_problem}
\min_{\boldsymbol p_{\overline{\omega^*}}} \, \frac{1}{2}\left[(\boldsymbol p^*_{\omega^*})^T \,\, (\boldsymbol p_{\overline{\omega^*}})^T\right]
\left[
\begin{array}{cc}
\boldsymbol S_{\omega^*, \omega^*} & \boldsymbol S_{\omega^*, \overline{\omega^*}}\\
\boldsymbol S_{\overline{\omega^*}, \omega^*} & \boldsymbol S_{\overline{\omega^*}, \overline{\omega^*}}
\end{array}
\right]
\left[\begin{array}{c}
\boldsymbol p^*_{\omega^*} \\
 \boldsymbol p_{\overline{\omega^*}}
\end{array}
\right] + (\boldsymbol a_{\omega^*}^T + \boldsymbol c^T \boldsymbol D_{\bullet, \omega^*}) \boldsymbol p^*_{\omega^*} +
(\boldsymbol a_{\overline{\omega^*}}^T + \boldsymbol c^T \boldsymbol D_{\bullet, \overline{\omega^*}}) \boldsymbol p_{\overline{\omega^*}}
\end{equation}
Assume for the sake of contradiction that the above claim does not hold. Since $\boldsymbol S_{\overline{\omega^*}, \overline{\omega^*}}$ is a principal submatrix of $\boldsymbol S$ and is hence positive definite, optimization problem (\ref{eq:restricted_problem}) has a unique minimizer. Let $\boldsymbol p'_{\overline{\omega^*}}$ be the minimizer of (\ref{eq:restricted_problem}). By the definition of $\omega^*$, we know that there exists $\lambda>0$ such that
$$ \left[
\begin{array}{c}
\boldsymbol p_{\omega^*} \\
\boldsymbol p_{\overline{\omega^*}} + \lambda (\boldsymbol p'_{\overline{\omega^*}} - \boldsymbol p^*_{\overline{\omega^*}})
\end{array}
\right] \, \in \, \Omega_{\kpd},$$
and clearly this vector also results in an objective value less than that of $p^*$. This is a contradiction. So the claim holds. And simple algebra leads to Equation (\ref{eq:partial_optimial}).
 \endproof
\begin{theorem}\label{th:profitbility_base}
Suppose assumptions (A1)-(A3) hold. Let $\boldsymbol p^*$ be an optimal solution of (\ref{eq:price_opt_const_base_linear_demand}), then $\boldsymbol p^*$ is profitable, i.e., $ \boldsymbol p^* \, \geq \, \boldsymbol c$.
\end{theorem}
\proof{}
Let $\omega^*$ be as defined in Lemma \ref{lm:partial optimality}. Since $\boldsymbol p^0 - \boldsymbol \delta \geq \boldsymbol c$ and $\boldsymbol \delta \geq 0$, we must have $\boldsymbol p^0 + \boldsymbol \delta \geq \boldsymbol p^0 \geq \boldsymbol c$. Therefore, for any $i \in \omega^*$, $p^*_i \geq c_i$, or equivalently $\boldsymbol p^*_{\omega^*} - \boldsymbol c_{\omega^*} \geq 0 $.  For $i \in \overline{\omega^*}$, we know from Lemma \ref{lm:partial optimality} that
\begin{eqnarray}
\boldsymbol p^{*}_{\overline {\omega^*}} - \boldsymbol c_{\overline{\omega^*}} & = & \boldsymbol S^{-1}_{\overline {\omega^*},\overline {\omega^*}}\left( \boldsymbol a_{\overline {\omega^*}}+\boldsymbol D_{\bullet, \overline {\omega^*}}^T \boldsymbol c -\boldsymbol S_{\overline {\omega^*},\omega^*} \boldsymbol p^*_{\omega^*} \right) - \boldsymbol c_{\overline{\omega^*}} \nonumber\\
& = & \boldsymbol S^{-1}_{\overline {\omega^*},\overline {\omega^*}}\left(\boldsymbol a_{\overline {\omega^*}}+ \boldsymbol D_{\bullet, \overline {\omega^*}}^T \boldsymbol c -\boldsymbol S_{\overline {\omega^*},\omega^*} \boldsymbol p^*_{\omega^*} - \boldsymbol S_{\overline {\omega^*},\overline {\omega^*}} \boldsymbol c_{\overline{\omega^*}}\right) \nonumber\\
& =& \boldsymbol S^{-1}_{\overline {\omega^*},\overline {\omega^*}}\left(\boldsymbol a_{\overline {\omega^*}}+ \boldsymbol D_{\omega^*, \overline {\omega^*}}^T \boldsymbol c_{\omega^*} + \boldsymbol D_{\overline{\omega^*},\overline {\omega^*}}^T \boldsymbol c_{\overline{\omega^*}}-(\boldsymbol D+\boldsymbol D^T)_{\overline {\omega^*}\omega^*} \boldsymbol p^*_{\omega^*} - (\boldsymbol D+\boldsymbol D^T)_{\overline {\omega^*},\overline {\omega^*}} \boldsymbol c_{\overline{\omega^*}}\right)\nonumber\\
& =& \boldsymbol S^{-1}_{\overline {\omega^*},\overline {\omega^*}}\left(\boldsymbol a_{\overline {\omega^*}}+ \boldsymbol D_{\omega^*, \overline {\omega^*}}^T \boldsymbol c_{\omega^*} -\boldsymbol D_{\overline {\omega^*},\omega^*} \boldsymbol p^*_{\omega^*} -\boldsymbol D_{\omega^*,\overline {\omega^*}}^T \boldsymbol p^*_{\omega^*} - \boldsymbol D_{\overline {\omega^*},\overline {\omega^*}} \boldsymbol c_{\overline{\omega^*}}\right)\nonumber\\
& =& \boldsymbol S^{-1}_{\overline {\omega^*},\overline {\omega^*}}\left(\boldsymbol a_{\overline {\omega^*}}+ \boldsymbol D_{\omega^*, \overline {\omega^*}}^T (\boldsymbol c_{\omega^*}- \boldsymbol p^*_{\omega^*}) -\boldsymbol D_{\overline {\omega^*},\omega^*} \boldsymbol p^*_{\omega^*} - \boldsymbol D_{\overline {\omega^*},\overline {\omega^*}} \boldsymbol c_{\overline{\omega^*}}\right)\nonumber\\
& =& \boldsymbol S^{-1}_{\overline {\omega^*},\overline {\omega^*}}\left(\boldsymbol a_{\overline {\omega^*}}+ \boldsymbol D_{\omega^*, \overline {\omega^*}}^T (\boldsymbol c_{\omega^*}- \boldsymbol p^*_{\omega^*}) +\boldsymbol D_{\overline {\omega^*},\omega^*} (\boldsymbol c_{\omega^*} - \boldsymbol p^*_{\omega^*} ) -\boldsymbol D_{\overline {\omega^*},\omega^*} \boldsymbol c_{\omega^*}- \boldsymbol D_{\overline {\omega^*},\overline {\omega^*}} \boldsymbol c_{\overline{\omega^*}}\right)\nonumber\\
& =& \boldsymbol S^{-1}_{\overline {\omega^*},\overline {\omega^*}}\left(\boldsymbol a_{\overline {\omega^*}} - \boldsymbol D_{\overline{\omega^*},\bullet} \boldsymbol c + \boldsymbol D_{\omega^*, \overline {\omega^*}}^T (\boldsymbol c_{\omega^*}- \boldsymbol p^*_{\omega^*}) +\boldsymbol D_{\overline {\omega^*},\omega^*} (\boldsymbol c_{\omega^*} - \boldsymbol p^*_{\omega^*} ) \right)\nonumber.
\end{eqnarray}
It is clear that $\boldsymbol S_{\overline {\omega^*},\overline {\omega^*}}$ is also an M-matrix and hence $\boldsymbol S^{-1}_{\overline {\omega^*},\overline {\omega^*}}$ is nonnegative (element-wise). By Assumption (A2), $\boldsymbol a-\boldsymbol D\boldsymbol c \geq \boldsymbol 0_n$, and hence $$(\boldsymbol a-\boldsymbol D\boldsymbol c)_{\overline{\omega^*}} \, = \, \boldsymbol a_{\overline {\omega^*}} - \boldsymbol D_{\overline{\omega^*},\bullet} \boldsymbol c \, \geq \, 0.$$ As we have shown above, $\boldsymbol p^*_{\omega^*} - \boldsymbol c_{\omega^*} \geq 0$. Also, by Assumption (A1) $\boldsymbol D_{\overline {\omega^*},\omega^*} \leq 0$ (element-wise) since it is an off-diagonal submatrix. Therefore,
$$ \boldsymbol D_{\overline {\omega^*},\omega^*}(\boldsymbol c_{\omega^*}-\boldsymbol p^*_{\omega^*}) \, \geq \, 0.$$
Similarly,
$$ \boldsymbol D^T_{\omega^*,\overline {\omega^*}}(\boldsymbol c_{\omega^*}-\boldsymbol p^*_{\omega^*}) \, \geq \, 0.$$
Overall, we have
$$\boldsymbol a_{\overline {\omega^*}} - \boldsymbol D_{\overline{\omega^*},\bullet} \boldsymbol c + \boldsymbol D_{\omega^*, \overline {\omega^*}}^T (\boldsymbol c_{\omega^*}- \boldsymbol p^*_{\omega^*}) +\boldsymbol D_{\overline {\omega^*},\omega^*} (\boldsymbol c_{\omega^*} - \boldsymbol p^*_{\omega^*} ) \, \geq \, 0.$$
And it follows that $\boldsymbol p^{*}_{\overline {\omega^*}} - \boldsymbol c_{\overline{\omega^*}} \geq 0$. The theorem thus holds readily.
 \endproof
Since $\boldsymbol S$ is symmetric positive definite, it has $n$ positive eigenvalues. Let its eignvalues be $\lambda_1\geq\lambda_2\geq\cdots\geq\lambda_n$. We collect some useful results regarding function $\mathcal Q(\boldsymbol p)$ in the following lemma.
\begin{lemma}\label{lm:linear_demand_properties}
Suppose that assumption (A1) holds. Then the following statements hold.
\begin{itemize}[leftmargin=.5in]
  \item[(a).] For any two price vectors $\boldsymbol p$ and $\boldsymbol p'$, and any $L \geq \lambda_1$
  $$\mathcal Q(\boldsymbol p') \, \leq \, \mathcal Q(\boldsymbol p) + \nabla \mathcal Q(\boldsymbol p)^T(\boldsymbol p'- \boldsymbol p) + \frac{1}{2} L \|\boldsymbol p- \boldsymbol p'\|_2^2. $$
  \item[(b).] For any price vector $\boldsymbol p$,
  $$\mathcal Q(\boldsymbol p) \, \leq \, \widehat{\mathcal Q} + \frac{1}{2 \lambda_n} \|\boldsymbol S (\boldsymbol p - \widehat{\boldsymbol p})\|_2^2. $$
  \item[(c).] The function $\mathcal Q(\boldsymbol p)$ is coercive, i.e.
  $$ \mathcal Q(\boldsymbol p ) \rightarrow +\infty \mbox{ as } \|\boldsymbol p\|_2 \rightarrow \infty.$$
\end{itemize}
\end{lemma}
\proof{}
Statement (a) can be shown using second order Taylor series expansion of a quadratic function at $\boldsymbol p$. In fact, we have:
$$ \mathcal Q(\boldsymbol p') \, = \, \mathcal Q(\boldsymbol p) + \nabla \mathcal Q(\boldsymbol p)^T(\boldsymbol p' - \boldsymbol p) + \frac{1}{2} (\boldsymbol p' - \boldsymbol p)^T H_{\mathcal Q}(\boldsymbol p) (\boldsymbol p' - \boldsymbol p),$$
where $H_{\mathcal Q}(\boldsymbol p)$ is the Hessian matrix at $\boldsymbol p$. It is clear that $H_{\mathcal Q}(\boldsymbol p) = \boldsymbol S$, and since $\boldsymbol S$ is symmetric positive definite with largest eigenvalue $\lambda_1$, we have
$$\mathcal Q(\boldsymbol p') \, \leq \, \mathcal Q(\boldsymbol p) + \nabla \mathcal Q(\boldsymbol p)^T(\boldsymbol p' - \boldsymbol p) + \frac{1}{2} \lambda_1\|\boldsymbol p' - \boldsymbol p\|_2^2 \, \leq \, \mathcal Q(\boldsymbol p) + \nabla \mathcal Q(\boldsymbol p)^T(\boldsymbol p' - \boldsymbol p) + \frac{1}{2} L\|\boldsymbol p' - \boldsymbol p\|_2^2 .$$
To show (b), we look at the Taylor series expansion at $\widehat{\boldsymbol p}$. Note that since $\widehat{\boldsymbol p}$ is the unconstrained minimizer of $\mathcal Q(\boldsymbol p)$, it follows that $\nabla \mathcal Q(\widehat{\boldsymbol p}) = \boldsymbol 0_n$. We have
\begin{eqnarray}
  \mathcal Q(\boldsymbol p) &=& \mathcal Q(\widehat{\boldsymbol p}) + \nabla \mathcal Q(\widehat{\boldsymbol p})^T(\boldsymbol p - \widehat{\boldsymbol p}) + \frac{1}{2} (\boldsymbol p - \widehat{\boldsymbol p})^T \boldsymbol S (\boldsymbol p - \widehat{\boldsymbol p}) \nonumber \\
  &=& \mathcal Q(\widehat{\boldsymbol p}) + \frac{1}{2} \left[(\boldsymbol p - \widehat{\boldsymbol p})^T\boldsymbol S\right]  \boldsymbol S^{-1} \left[ \boldsymbol S (\boldsymbol p - \widehat{\boldsymbol p})\right] \nonumber\\
  & \leq & \widehat{\mathcal Q} + \frac{1}{2\lambda_n} \|\boldsymbol S(\boldsymbol p - \widehat{\boldsymbol p})\|_2^2,\nonumber
\end{eqnarray}
where the inequality is due to the fact that the eigenvalues of $\boldsymbol S^{-1}$ are $\frac{1}{\lambda_n} \geq \frac{1}{\lambda_{n-1}}\geq \cdots\geq \frac{1}{\lambda_1}$.\\
Statement (c) is a direct consequence of the positive definiteness of $\boldsymbol S$.
 \endproof
Note that the coercivity of $\mathcal Q(\boldsymbol p)$ implies that it has bounded level sets, i.e., for any given $c\in \mathbb R$, the level set $L_c \triangleq \{\boldsymbol p \in \mathbb R^n \, |\, \mathcal Q(\boldsymbol p) \leq c\}$ is bounded (and closed).
\subsection{Combinatorial Nature}
It is clear that $\Omega_{\kpd}$ is not convex. On the other hand, we can show that it is in fact the union of finitely  many polyhedra. To see this, let a triplet of index sets $(\alpha,\beta,\gamma)$ be a partition of $\{1,2,\cdots,n\}$, satisfying:
\begin{equation} \label{eq:paritition_condition}
|\beta|+ |\gamma| \leq k,\,\,\alpha \cup \beta \cup \gamma = \{1,2,\cdots,n\}, \,\,\alpha \cap \beta = \emptyset, \,\,\beta \cap \gamma = \emptyset, \mbox{ and } \alpha \cap \gamma = \emptyset.
\end{equation}
Let the family of all partitions satisfying (\ref{eq:paritition_condition}) be denoted by $\Im$. Let the polyhedron $\mathcal F_{\alpha,\beta,\gamma}$ be defined as:
\begin{equation}\label{eq:partition_polyhedron}
\mathcal F_{\alpha,\beta,\gamma} \, \triangleq \, \left\{
\boldsymbol p \in \mathbb R^n \, \left|\,
\begin{array}{ll}
p_i \, = \, p_i^0, & \forall \, i \in \alpha,\\
p_i \, \geq  \, p_i^0+\delta_i, & \forall \, i \in \beta, \\
p_i \, \leq \, p_i^0-\delta_i, & \forall \, i \in \gamma.
\end{array}
 \right.\right\}
\end{equation}
It is easy to verify that $$\Omega_{\kpd} = \bigcup_{(\alpha,\beta,\gamma) \in \Im} \mathcal F_{\alpha,\beta,\gamma}.$$ Throughout the rest of the paper, we make an additional assumption:
\begin{itemize}[leftmargin=.5in]
  \item[(A4).] $\delta_i$'s are all positive, i.e., $\delta_i > 0$ for all $i=1,\cdots,n$.
\end{itemize}
Assumption (A4) is not restrictive in practice. In fact, in reality any price change is at least by 1 cent, therefore (A4) is valid. Under Assumption (A4), $\mathcal F_{\alpha,\beta,\gamma}$'s are non-overlapping. We can define an optimization problem over each polyhedron $\mathcal F_{\alpha,\beta,\gamma}$ as follows.
\begin{equation}\label{eq:restricted_opt_anypartition}
\begin{array}{rll}
\overline{\boldsymbol p}_{\alpha,\beta,\gamma} \, \triangleq \, \displaystyle{ \argmin_{\boldsymbol p} } & \mathcal Q(\boldsymbol p) \\[5pt]
\mbox{s.t.} & p_i \, = \, p_i^0, & \forall \, i \in \alpha, \\
&p_i \, \geq  \, p_i^0+\delta_i, & \forall \, i \in \beta, \\
&p_i \, \leq \, p_i^0-\delta_i, & \forall \, i \in \gamma.
\end{array}
\end{equation}
Due to Assumption (A1), optimization problem (\ref{eq:restricted_opt_anypartition}) has a unique solution, denoted by $\overline{\boldsymbol p}_{\alpha,\beta,\gamma}$. With Assumption (A4), it is easy to verify that $\overline{\boldsymbol p}_{\alpha,\beta,\gamma}$ is a local optimal solution of (\ref{eq:price_opt_const_base_linear_demand}). Therefore, a global optimal solution of (\ref{eq:price_opt_const_base_linear_demand}) can be found by enumerating all possible partitions satisfying (\ref{eq:paritition_condition}). However, the number of partitions in $\Im$ is exponentially large with respect to the number of product $n$. In fact, the number of possible partitions in $\Im$ is
$$\sum_{\imath=1}^k  { {n}\choose{\imath}} 2^{\imath}. $$ It is hence not practical to enumerate all the possible partitions. This observation demonstrates the combinational nature of the price optimization problem (\ref{eq:price_opt_const_base_linear_demand}). Therefore, in the following section we study an algorithm that finds a good near-optimal solution.

\section{A Gradient Projection Algorithm for Price Optimization} \label{sec:algorithm}
In this section, we study a gradient projection method for the price optimization problem (\ref{eq:price_opt_const_base_linear_demand}). The gradient projection method, due to its simplicity, is a widely used optimization algorithm. It can be regarded as a natural extension of the classical gradient descent method for unconstrained optimization problems to the constrained ones. Typically, for provable convergence to an optimal solution, the optimization problem needs to have both a convex objective function and a convex feasible set. In problem (\ref{eq:price_opt_const_base_linear_demand}), the objective function is convex but the feasible set is not. On the other hand, as we will show, the feasible set, albeit nonconvex, allows for an easily computable (in polynomial time) algorithm to compute a projection onto it. In this research, we will show both theoretically and empirically that a gradient projection method can perform very well for problem (\ref{eq:price_opt_const_base_linear_demand}) even though a provable optimal solution is not guaranteed to be found. Our method is motivated by \cite{bertsimas2016best}, where the authors studied regression problems with zero norm constraints. To present our algorithm, we start with some results regarding projections onto the feasible set $\Omega_{\kpd}$.
\subsection{Projection Onto the Feasible Set}
We first consider the one-dimensional projection problem onto $\pdi$. Given any $q_i \in \mathbb R$, we aim to project it onto $\mathcal P^i_{\delta_i}$, or equivalently find a $ \widetilde p_i$ to
\begin{equation}\label{eq:1d_projection}
\min_{p_i \in \pdi} (p_i -q_i)^2.
\end{equation}
Since $\pdi$ is not convex, the optimization problem (\ref{eq:1d_projection}) does not necessarily have a unique solution, i.e., the projection is not necessarily a singleton. Let the set of solutions to (\ref{eq:1d_projection}) be denoted as $\Pi_{\pdi}(q_i)$. It is clear that
\begin{equation}\label{eq:1d_projection_closedform}
\argmin_{p_i \in \pdi} \,(p_i -q_i)^2 \, \triangleq \, \Pi_{\mathcal P^i_{\delta_i}}(q_i) \, = \,
\left\{
\begin{array}{ll}
\{q_i\} & \mbox{if } q_i \in \pdi\\
\{p_i^0 + \delta_i\} & \mbox{if } p^0_i + \frac{1}{2}\delta_i < q_i < p^0_i + \delta_i \\
\{p_i^0, p_i^0 + \delta_i\} & \mbox{if } q_i = p^0_i + \frac{1}{2}\delta_i \\
\{p_i^0\} & \mbox{if } p_i^0-\frac{1}{2}\delta_i < q_i < p_i^0+\frac{1}{2}\delta_i\\
\{p_i^0, p_i^0-\delta_i\} & \mbox{if } q_i = p^0_i - \frac{1}{2}\delta_i\\
\{p_i^0 - \delta_i\} & \mbox{if } p^0_i - \delta_i < q_i < p^0_i - \frac{1}{2}\delta_i
\end{array}
\right.
\end{equation}
We can see that for any $\boldsymbol q \in \mathbb R^n$, the projection of it onto $\mathcal P_{\boldsymbol \delta} \subseteq \mathbb R^n$ is
$$ \Pi_{\mathcal P_{\boldsymbol \delta}}(\boldsymbol q) = \bigotimes_{i=1}^n \Pi_{\pdi} (q_i).$$
Since $\Pi_{\mathcal P^i_{\delta_i}}(q_i)$ is not necessarily a singleton, we define a set-valued map
$$ \phi_i(q_i): q_i \mapsto \Pi_{\pdi}(q_i).$$
The graph of $\phi_i$ is defined as
$$\mbox{gph } \phi_i \, \triangleq \, \{ (q_i, g_i) \in \mathbb R^2 | g_i \in \phi_i(q_i)\}.$$
Recall that a set-valued map $\phi : \mathbb R^n \mapsto \mathbb R^n$ is closed on $\mathbb R^n$ if for any point $\boldsymbol x \in \mathbb R^n$, the following implication holds:
\begin{equation} \label{eq:closed_set_valued_map}
\left.
\begin{array}{r}
\lim_{j\rightarrow \infty}\boldsymbol x^j = \boldsymbol x \\[5pt]
\lim_{j\rightarrow \infty}\boldsymbol y^j = \boldsymbol y \\[5pt]
\boldsymbol y^j \in \phi(\boldsymbol x^j), \,\, \forall\, j
\end{array}
\right\} \, \Rightarrow \, \boldsymbol y \in \phi(\boldsymbol x).
\end{equation}
It is well known that a set valued map $\phi$ is closed on $\mathbb R^n$ if and only if its graph $\mbox{gph } \phi$ is a closed set in $\mathbb R^n$. We have the following lemma regarding the closedness of $\phi_i$.
\begin{lemma}\label{lm:closed_single_projection}
	For all $i=1,\cdots,n$, $\phi_i$ is a closed map.
\end{lemma}
\proof{}
	The graph of $\phi_i$ on the $q_i$-$g_i$ plane is the union of 2 rays and 3 line segments as follows:
	\begin{eqnarray}
	\mbox{gph } \phi_i  & = & \{(q_i,g_i) | q_i = g_i, q_i \geq p_i^0+\delta_i \} \cup \{(q_i,g_i) | q_i = g_i, q_i \leq p_i^0-\delta_i \} \nonumber \\
	&& \cup \{(q_i,g_i) | g_i = p_i^0+\delta_i, p_i^0+\frac{1}{2}\delta_i \leq q_i \leq p_i^0+\delta_i \} \nonumber\\
	&& \cup \{(q_i,g_i) | g_i = p_i^0, p_i^0-\frac{1}{2}\delta_i \leq q_i \leq p_i^0+\frac{1}{2}\delta_i \} \nonumber\\
	&& \cup \{(q_i,g_i) | g_i = p_i^0-\delta_i, p_i^0-\delta_i \leq q_i \leq p_i^0-\frac{1}{2}\delta_i \}. \nonumber
	\end{eqnarray}
	Since rays and line segments are closed, the union of these five sets is also closed. This concludes the proof.
 \endproof
\begin{figure}
	\centering
	\begin{subfigure}{.5\textwidth}
		\centering
		\includegraphics[width=0.8\linewidth]{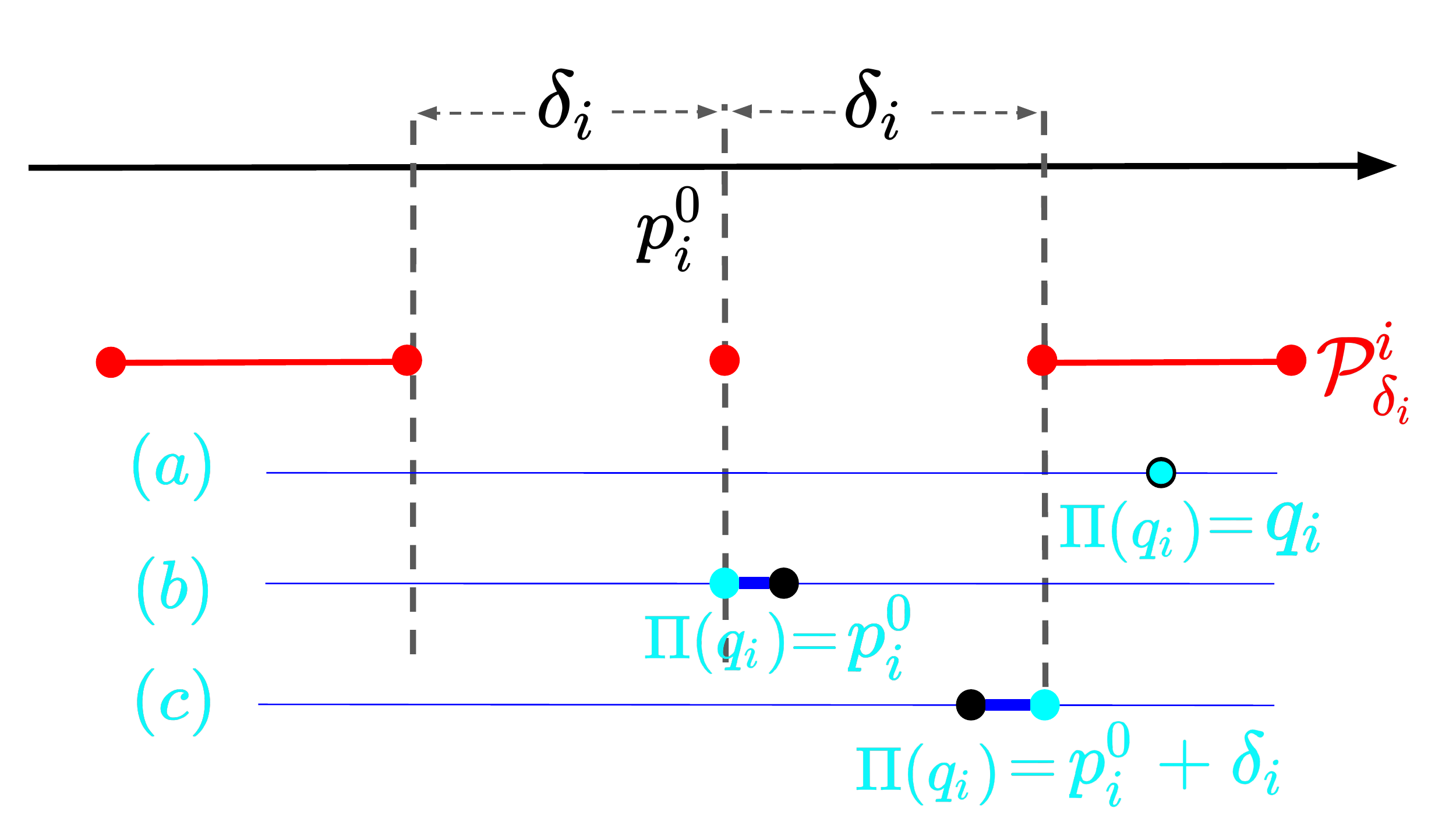}
		\caption{ \centering Projection on to $\pdi$.}
		\label{fig:one_dimensional_projection}
	\end{subfigure}%
	\begin{subfigure}{.5\textwidth}
		\centering
		\includegraphics[width=1.0\linewidth]{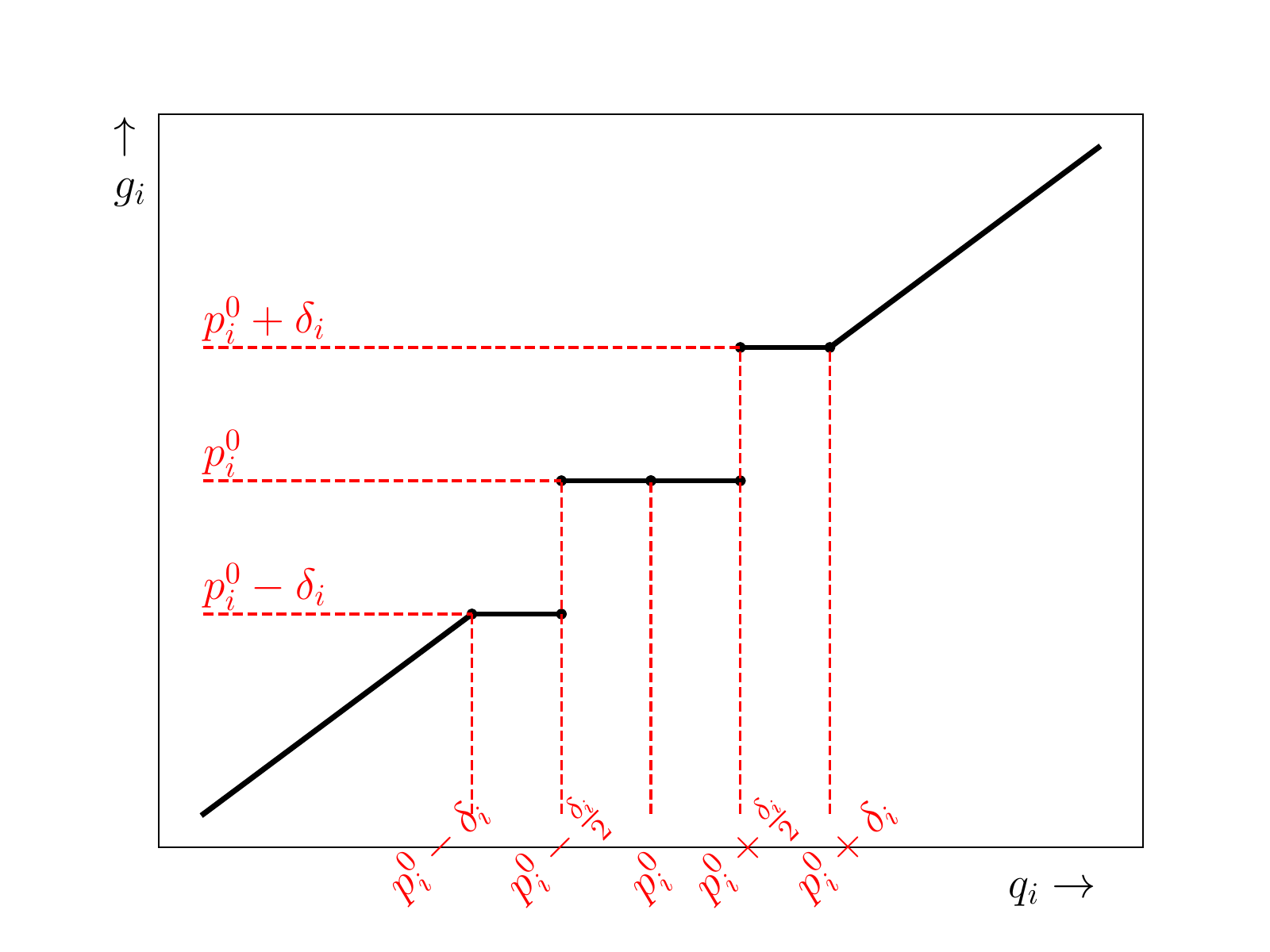}
		\caption{ \centering Graph of set-valued map $\phi_i$.}
		\label{fig:projection_map}
	\end{subfigure}
	\caption{\centering Projection onto the feasible set.}
\end{figure}

Define $d(q_i,\pdi) \, \triangleq \, \displaystyle{ \min_{p_i \in \pdi} } (p_i - q_i)^2,$
i.e., the optimal value of (\ref{eq:1d_projection}) or the squared distance from $q_i$ to set $\pdi$. It is straightforward to verify that
\begin{equation}\label{eq:distance_upperbound}
d(q_i,\pdi) \, \leq \, \frac{\delta_i^2}{4}.
\end{equation}
In fact, we have
\begin{equation}\label{eq:distance_closed_form}
d(q_i,\pdi) \, = \, \left\{
\begin{array}{ll}
0 & \text{if } q_i \geq p^0_i + \delta_i \\
(p^0_i + \delta_i - q_i)^2 & \text{if } p^0_i+\frac{1}{2}\delta_i \leq  q_i \leq p^0_i + \delta_i \\
(q_i-p^0_i)^2 & \text{if } p^0_i \leq  q_i \leq p^0_i +\frac{1}{2} \delta_i \\
(p^0_i-q_i)^2 & \text{if } p^0_i-\frac{1}{2} \delta_i \leq q_i \leq p^0_i \\
(q_i - p^0_i + \delta_i)^2 & \text{if } p^0_i- \delta_i \leq q_i \leq  p^0_i-\frac{1}{2}\delta_i\\
0 & \text{if } q_i \leq p^0_i - \delta_i
\end{array}
\right.
\end{equation}
It is clear that $d(q_i,\pdi)$, when regarded as a function of $q_i$ is continuous. For each $i \in \{1,\cdots,n \}$, define:
\begin{equation}\label{eq:Deltas}
\Delta_i \, \triangleq \, \left(p^0_i - q_i\right)^2-d(q_i,\pdi).
\end{equation}
By the definition of $d(q_i,\pdi)$, and the fact that $p^0_i \in \pdi$, it is clear that $\Delta_i \geq 0$. We have the following lemma.
\begin{lemma}\label{lm:zero_delta}
For each $i=1,\cdots,n$, Let $\Delta_i$ be defined by (\ref{eq:Deltas}), it holds that $\Delta_i = 0$ if and only if $ p_i^0-\frac{1}{2} \delta_i \leq q_i \leq p_i^0+\frac{1}{2} \delta_i$.
\end{lemma}
\proof{}
It is clear that $\Delta_i = 0$ if and only if $p_i^0 \in \Pi_{\pdi} (q_i)$. From Equation (\ref{eq:1d_projection_closedform}), it is clear that $p_i^0 \in \Pi_{\pdi} (q_i)$ if and only if $p_i^0-\frac{1}{2} \delta_i \leq q_i \leq p_i^0+\frac{1}{2} \delta_i$.
 \endproof
We next study the projection of a vector $\boldsymbol q\in \mathbb R^n$ onto the set $\Omega_{\kpd}$, i.e., the following optimization problem:
	\begin{equation}\label{eq:sub-prob}
    \begin{array}{rl}
	\displaystyle{ \min_{ \boldsymbol p } } & || \boldsymbol p - \boldsymbol q||_2^2  \\ [5pt]
    \mbox{s.t.} & ||\boldsymbol p- \boldsymbol p^0 ||_0 \, \leq \, k \\
     & \boldsymbol p \in \mathcal P_{\boldsymbol \delta}
    \end{array}
	\end{equation}
\begin{proposition}\label{prop:sub_problem_prop}
Assume $\Delta_{(1)} \geq \Delta_{(2)} \geq \cdots \geq \Delta_{(n)}$, i.e., the ordered values of $\Delta_i$'s defined in (\ref{eq:Deltas}), let $\widetilde{\boldsymbol p}$ satisfies
	\begin{equation}\label{eq:sub_prob_solution}
	\widetilde p_i \, \begin{cases} \in \Pi_{\pdi} (q_i), & \mbox{if } i \in \{ (1),...,(k)\}, \\
	= p_i^0, & \mbox{otherwise}.  \end{cases}
	\end{equation}
then $\widetilde{\boldsymbol p}$ is an optimal solution to (\ref{eq:sub-prob}).
\end{proposition}
\proof{} It is straightforward to verify that $\widetilde{\boldsymbol p}$ given by (\ref{eq:sub_prob_solution}) is a feasible solution to (\ref{eq:sub-prob}). Let $\mathcal S_k$ be the collection of all subsets of $\{1,\cdots,n\}$ with cardinality less than or equal to $k$. For each $\sigma \in \mathcal S_k$, we define the following optimization problem
\begin{equation}\label{eq:subset_S}
	\begin{array}{rl}
	V(\sigma) \, \triangleq \, \displaystyle{ \min_{ \boldsymbol p } }& \displaystyle { \sum_{i \notin \sigma} (q_i - p_i^0)^2 + \sum_{i \in \sigma} (q_i - p_i)^2  } \\ [5pt]
	\mbox{s.t.} &  p_i = p^0_i, \,\, \forall \, i \notin \sigma,\\
    & \boldsymbol p \in \mathcal P_{\boldsymbol \delta}.
	\end{array}
\end{equation}
Let $p(\sigma)$ be given by
\begin{equation}\label{eq:projection}
p(\sigma)_i \,  \left\{
\begin{array}{ll}
\in \Pi_{\mathcal P_{\delta_i}} (q_i) & \mbox{ if } i \in \sigma \\
= \, p^0_i & \mbox{ if } i\notin \sigma .
\end{array}
\right.
\end{equation}
It is clear that $p(\sigma)$ is a feasible solution for (\ref{eq:subset_S}). Moreover, given any vector $\boldsymbol p$ feasible for  (\ref{eq:subset_S}), we must have
$$\sum_{i \notin \sigma} (q_i - p_i^0)^2 + \sum_{i \in \sigma} (q_i - p_i)^2   \, \geq \, \sum_{i \notin \sigma} (q_i - p_i^0)^2 + \sum_{i \in \sigma} d(q_i, \pdi)  .$$
Therefore, $p(\sigma)$ is an optimal solution of (\ref{eq:subset_S}), with the optimal value
\begin{eqnarray}
V^*(\sigma) & = & \sum_{i \notin \sigma} (q_i - p_i^0)^2 + \sum_{i \in \sigma} d(q_i, \pdi) \nonumber\\
& = & \sum_{i =1}^n (q_i - p_i^0)^2 + \sum_{i \in \sigma} \left[d(q_i, \pdi) - (q_i - p_i^0)^2\right] \nonumber\\
& = & \|\boldsymbol q - \boldsymbol p^0\|_2^2  - \sum_{i \in \sigma} \Delta_i.\nonumber
\end{eqnarray}
Notice that (\ref{eq:sub-prob}) is equivalent to
$$\min_{\sigma\in \mathcal S_k} V^*(\sigma). $$
Since $\Delta_{(1)} \geq \Delta_{(2)} \geq \cdots \geq \Delta_{(n)}$, for any $\sigma \in \mathcal S_k$
$$ V^*(\sigma) \, \geq \, \|\boldsymbol q - \boldsymbol p^0\|^2 - \sum_{i \in \{(1),\cdots,(k)\}} \Delta_i,$$
with the equality holding when $\sigma$ is taken to be $\{(1),\cdots,(k)\}$. This implies that $\widetilde{\boldsymbol p}$ given in (\ref{eq:sub_prob_solution}) is an optimal solution of (\ref{eq:sub-prob}).
 \endproof
The above proposition provide an efficient algorithm to find a projection onto $\Omega_{\kpd}$ for any given $\boldsymbol q\in \mathbb R^n$. Due to the non-convexity of set $\Omega_{\kpd}$, the solution set of (\ref{eq:sub-prob}) is not necessarily a singleton. For example, there can be ties among $\Delta_i$'s, leading to multiple index sets of $k$ largest $\Delta_i$'s and hence multiple optimal solutions. In addition, the projection $\Pi_{\pdi}(q_i)$ itself is not necessarily a singleton. On the other hand, based on the proof of Proposition \ref{prop:sub_problem_prop} it is straightforward to prove the next lemma, which provides a property of a solution to (\ref{eq:sub-prob}).
\begin{lemma}\label{lm:projection_omega_construction_necessary}
Let $\widetilde{\boldsymbol p} \in \mathbb R^n$ define $\sigma \triangleq \{i \, |\, \widetilde{p}_i \neq p^0_i\}$. The following statements hold.
\begin{itemize}
    \item[(i)] If $|\sigma|= k$, then $\widetilde{\boldsymbol p}$ is an optimal solution to (\ref{eq:sub-prob}) if and only if the following conditions hold:
$$
\left\{
\begin{array}{ll}
\widetilde p_i \in \Pi_{\pdi} (q_i),& \forall \, i \in \sigma,\\
\Delta_i \geq \Delta_j, & \forall \, i\in \sigma, j\notin \sigma.
\end{array}
\right.
$$
\item[(ii)] If $|\sigma|< k$, then $\widetilde{\boldsymbol p}$ is an optimal solution to (\ref{eq:sub-prob}) if and only if the following conditions hold:
$$
\left\{
\begin{array}{ll}
\widetilde p_i \in \Pi_{\pdi} (q_i),& \forall \, i \in \sigma,\\
\Delta_i = 0, & \forall \, i \notin \sigma.
\end{array}
\right.
$$
\end{itemize}
\end{lemma}

In the rest of the paper, we denote the solution set to problem \eqref{eq:sub-prob} by the notation $\mathcal H_{\kpd}(\boldsymbol q)$. We next study the set-valued map
$$\varphi(\boldsymbol q): \boldsymbol q \mapsto \mathcal H_{\kpd}(\boldsymbol q).$$
We have the following result regarding $\varphi(\boldsymbol q)$.
\begin{lemma}\label{lm:closed_map}
	$\varphi(\boldsymbol q)$ is a closed set-valued map on $\mathbb R^n$.
\end{lemma}
\proof{}
	We prove the lemma by verifying the definition. Take an arbitrary $\boldsymbol q \in \mathbb R^n$, let sequences $\left\{\boldsymbol q^j\right\}_{j=1,2,\cdots}$ and $\left\{\boldsymbol g^j\right\}_{j=1,2,\cdots}$ satisfy $$\lim_{j\rightarrow \infty}\boldsymbol q^j = \boldsymbol q, \lim_{j\rightarrow \infty}\boldsymbol g^j = \boldsymbol g, \mbox{ and } \boldsymbol g^j \in \mathcal H_{\kpd}(\boldsymbol q^j) \,\, \forall\, j.$$ Define $\Delta^j_i \triangleq (p_i^0 - q^j_i)^2 - d(q_i^j, \pdi)$. For each $j=1,2,\cdots$, let $\sigma^j \triangleq \{i \, | \, g^j_i \neq p^0_i\}$, it is clear that $|\delta| \leq k$. By Lemma \ref{lm:projection_omega_construction_necessary} it holds that
	$$ g^j_i \in \pdi(q^j_i), \forall \, i \in \sigma^j, \,\, \mbox{ and } \Delta^j_i \geq \Delta^j_{i'}, \forall\,\, i\in \sigma^j, i'\notin \sigma^j.$$
	Since there are only finitely many possible subsets of $\{1,2,\cdots,n\}$, at least one subset appears infinitely mainly times in the sequence $\{\sigma^j\}_{j=1,2,\cdots}$. Therefore, by working with a subsequence if necessary, we can assume without loss of generality that $\sigma^j = \sigma^*$ for all $j=1,2,\cdots$. It is clear that for any $i\in \sigma^*$ and $i'\notin \sigma^*$, by the continuity of $d(q_i, \pdi)$, we have
	\begin{eqnarray}
	\Delta_i & = & (p_i^0 - q_i)^2 - d(q_i, \pdi) \, = \, \lim_{j\rightarrow \infty} \left[(p_i^0 - q_i^j)^2 - d(q_i^j, \pdi)\right] \nonumber \\
	& \geq & \lim_{j\rightarrow \infty} \left[(p_{i'}^0 - q_{i'}^j)^2 - d(q_{i'}^j, \mathcal P^{i'}_{\delta_{i'}}) \right]  \, = \, (p_{i'}^0 - q_{i'})^2 - d(q_{i'}, \mathcal P^{i'}_{\delta_{i'}})\, = \, \Delta_{i'} \nonumber
	\end{eqnarray}
	In addition, in case $|\sigma^*| < k$, we have:
	$$ \Delta_i \,  = \,  (p_i^0 - q_i)^2 - d(q_i, \pdi) \, = \, \lim_{j\rightarrow \infty} \left[(p_i^0 - q_i^j)^2 - d(q_i^j, \pdi) \right] \, = \, 0, \,\, \forall \, i \notin \sigma^*. $$
	Moreover, by Lemma \ref{lm:closed_single_projection}, i.e., the closedness of $\phi_i$, we have
	$$ g_i \in \pdi(q_i), \forall \, i \in \sigma^*. $$
	It is also clear that $g_{i'} = p_{i'}^0, \forall \, i' \notin \sigma^*$. Therefore, by Lemma \ref{lm:projection_omega_construction_necessary}, we have $\boldsymbol g \in \mathcal H_{\kpd}(\boldsymbol q)$. This concludes the proof.
 \endproof
\begin{remark}
Proposition \ref{prop:sub_problem_prop} is a generalization of Proposition 3 of \cite{bertsimas2016best}, that deals with the special case where $\boldsymbol \delta = \boldsymbol 0_{n}$  and $\boldsymbol p^0 = \boldsymbol 0_n$.
\end{remark}
\begin{remark}
When upper and lower bounds on individual products exist, the projection of any $q_i \in \mathbb R$ on to set $\mathcal P^i_{\delta_i,l_i,u_i}$ is given by
$$
\Pi_{\mathcal P^i_{\delta_i,l_i,u_i}}(q_i) \, = \,
\left\{
\begin{array}{ll}
\{q_i\} & \mbox{if } q_i \in \mathcal P_{\delta_i,l_i,u_i}\\
\{l_i\} & \mbox{if } q_i < l_i \\
\{p_i^0 - \delta_i\} & \mbox{if } p^0_i-\delta_i < q_i < p^0_i-\frac{1}{2} \delta_i\\
\{p_i^0 - \delta_i, p_i^0\} & \mbox{if } q_i = p^0_i-\frac{1}{2} \delta_i\\
\{p_i^0\} & \mbox{if } p^0_i-\frac{1}{2} \delta_i < q_i < p^0_i+\frac{1}{2} \delta_i\\
\{p_i^0, p_i^0 +\delta_i\} & \mbox{if } q_i = p^0_i+\frac{1}{2} \delta_i\\
\{p_i^0 + \delta_i\} & \mbox{if } p^0_i+\frac{1}{2} \delta_i \leq q_i \leq p^0_i+\delta_i\\
\{u_i\} & \mbox{if } q_i > u_i
\end{array}
\right.
$$
Proposition \ref{prop:sub_problem_prop} as well as Lemma \ref{lm:closed_map} can then be extended to include this case. The proofs are straightforward and hence omitted.
\end{remark}
\subsection{A Gradient Projection Algorithm}
The gradient projection algorithm can be regarded as minimizing an upper approximation at current solution of $\mathcal Q(\boldsymbol p)$ given by Lemma \ref{lm:linear_demand_properties} item (a). In particular, at iteration $t$, let the current solution be $\boldsymbol p^t$, for an $L>\lambda_1$, we define
\begin{equation}\label{eq:Q_upper}
\widetilde{\mathcal Q}_L(\boldsymbol p;\boldsymbol p^t) \, \triangleq \, \mathcal Q(\boldsymbol p^t) + \nabla \mathcal Q(\boldsymbol p^t)^T(\boldsymbol p- \boldsymbol p^t) + \frac{1}{2} L\|\boldsymbol p- \boldsymbol p^t\|_2^2.
\end{equation}
It is easy to see that
\begin{eqnarray}
\argmin_{\boldsymbol p\in \Omega_\kpd} \widetilde{ \mathcal Q }_L(\boldsymbol p; \boldsymbol p^t)
	 & = & \argmin_{\boldsymbol p \in \Omega_\kpd} \left( \frac{L}{2} \left\|\boldsymbol p - \left(\boldsymbol \lambda_1 -\frac{1}{L} \nabla \mathcal Q(\boldsymbol p^t)\right) \right \|_2^2 - \frac{1}{2L} \|\nabla \mathcal Q(\boldsymbol p^t) \|_2^2 + \mathcal Q(\boldsymbol p^t) \right) \nonumber \\
	 &  = &	\argmin_{\boldsymbol p \in \Omega_\kpd} \left\|\boldsymbol p - \left(\boldsymbol p^t- \frac{1}{L} \nabla \mathcal Q(\boldsymbol p^t)\right)  \right\|_2^2 \nonumber \\
	 & = & \mathcal H_{k,\boldsymbol \delta, \boldsymbol p^0}\left(\boldsymbol p^t- \frac{1}{L} \nabla \mathcal Q(\boldsymbol p^t)\right),\label{eq:proximal}
\end{eqnarray}
i.e., the projection of $\boldsymbol p^t- \frac{1}{L} \nabla \mathcal Q(\boldsymbol p^t)$ onto the set $\Omega_{\kpd}$.
The algorithm is given below.\\

\begin{center}
\noindent \fbox{	
\begin{minipage}{0.8\textwidth}
\textbf{Gradient Projection Algorithm (GPA)} \\
Input: $\mathcal Q(\boldsymbol p)$, parameter: $L>\lambda_1$ and convergence tolerance $\varepsilon$.\\
Output: A first-order stationary solution $\boldsymbol p^\infty$.
\begin{itemize}[leftmargin=.5in]
\item[1.] Initialize with $\boldsymbol p^1$ such that $||\boldsymbol p^1-\boldsymbol p^0||_0 \leq k$.
\item[2.] For $t \geq 1 $, apply \eqref{eq:sub_prob_solution} to obtain $\boldsymbol p^{t+1}$, i.e.,
\begin{equation}\label{eq:iter}
\boldsymbol p^{t+1} \in \mathcal H_{\kpd}\left(\boldsymbol p^{t}- \frac{1}{L} \nabla \mathcal Q(\boldsymbol p^{t})\right).
\end{equation}
\item[3.] Repeat step 2, until $\mathcal Q(\boldsymbol p^t)- \mathcal Q(\boldsymbol p^{t+1}) \leq \varepsilon$.
\end{itemize}
\end{minipage}
}
\end{center}
\subsection{Analysis of the Gradient Projection Algorithm}
In this subsection, we analyze the GPA. We aim to answer the following questions:
\begin{itemize}[leftmargin=.5in]
  \item[1.] What kind of solutions is the GPA designed to find? How good are these solutions in terms of the objective values?
  \item[2.] Does the GPA indeed converge to a solution it is designed to find?
\end{itemize}
It is clear that the GPA essentially aims to find a fixed point of the (set-valued) map
$$\boldsymbol p \, \mapsto \, \mathcal H_{\kpd}\left(\boldsymbol p - \frac{1}{L} \nabla \mathcal Q(\boldsymbol p)\right).$$
We first study how such a fixed point is related to an optimal solution of problem (\ref{eq:price_opt_const_base_linear_demand}). A technical lemma regarding a property of $\mathcal H_{\kpd}(\cdot)$ is stated below.
\begin{lemma}\label{lm:H_property}
Let $\boldsymbol p' \in \mathcal H_{\kpd}\left(\boldsymbol p - \frac{1}{L} \nabla \mathcal Q(\boldsymbol p)\right)$, then
\begin{equation}\label{eq:decreasing}
\mathcal Q( \boldsymbol p) - \mathcal Q(\boldsymbol p') \, \geq \, \frac{L-\lambda_1}{2} \|\boldsymbol p' - \boldsymbol p\|_2^2.
\end{equation}
\end{lemma}
\proof{}
As we can see from equation (\ref{eq:proximal}),
$$ \boldsymbol p' \, \in \, \argmin_{\boldsymbol \eta \in \Omega_{\kpd}} \widetilde{\mathcal Q}_L(\boldsymbol \eta;\boldsymbol p).$$
Therefore
\begin{eqnarray}
\mathcal Q(\boldsymbol p) \, = \,  \widetilde{\mathcal Q}_L(\boldsymbol p;\boldsymbol p) & \geq & \widetilde{\mathcal Q}_L(\boldsymbol p';\boldsymbol p)\nonumber\\
 & = & \mathcal Q(\boldsymbol p) + (\nabla \mathcal Q(\boldsymbol p))^T(\boldsymbol p' - \boldsymbol p) + \frac{L}{2}\|\boldsymbol p' - \boldsymbol p\|_2^2 \nonumber\\
& = & \mathcal Q(\boldsymbol p) + (\nabla \mathcal Q(\boldsymbol p))^T(\boldsymbol p' - \boldsymbol p)+\frac{\lambda_1}{2}\|\boldsymbol p' - \boldsymbol p\|_2^2 + \frac{L-\lambda_1}{2}\|\boldsymbol p' - \boldsymbol p\|_2^2 \nonumber\\
&\geq & \mathcal Q(\boldsymbol p') + \frac{L-\lambda_1}{2}\|\boldsymbol p' - \boldsymbol p\|_2^2, \nonumber
\end{eqnarray}
where the second inequality is due to Lemma \ref{lm:linear_demand_properties} item (a). Thus, we have
\begin{equation}
\mathcal Q(\boldsymbol p)  - \mathcal Q(\boldsymbol p') \, \geq  \, \frac{L- \lambda_1}{2} ||\boldsymbol p'- \boldsymbol p||_2^2.
\end{equation}
 \endproof
We define a first-order stationary point as follows.
\begin{definition}\label{def:first_order_stationary}
[First-Order Stationary Point] Given an $L > \lambda_1$, a vector $\boldsymbol p^{\infty} \in \mathbb{R}^n$ is said to be a first-order
stationary point of problem \eqref{eq:price_opt_const_base_linear_demand}  if it satisfies the following fixed point equation:
\begin{equation}
\begin{aligned}
\boldsymbol p^\infty \, \in \, \mathcal H_{\kpd}\left(\boldsymbol p^\infty- \frac{1}{L} \nabla \mathcal Q(\boldsymbol p^\infty)\right).
\end{aligned}	
\label{eq:stationary}
\end{equation}
\end{definition}
\noindent Note that (\ref{eq:stationary}) implies that $\boldsymbol p^\infty \in \Omega_{\kpd}$ by the definition of $\mathcal H_{k,\boldsymbol \delta, \boldsymbol p^0}(\cdot)$. Let
\begin{equation}\label{eq:index_set_sationary}
\alpha^\infty \triangleq \{i \, | \, p^\infty_i = p^0_i\},\,\, \beta^\infty \triangleq \{i \, | \, p^\infty_i \geq p^0_i+\delta_i\}, \,\, \mbox{and}\,\, \gamma^\infty \triangleq \{i \, | \, p^\infty_i \leq p^0_i - \delta_i\}.
\end{equation}
Due to assumption (A4), the index sets $\alpha^\infty$, $\beta^\infty$, and $\gamma^\infty$ are uniquely defined given $\boldsymbol p^\infty$. With the index sets, define the following optimization problem
\begin{equation}\label{eq:restricted_opt}
\begin{array}{rlll}
\displaystyle{ \min_{\boldsymbol p} } & \mathcal Q(\boldsymbol p) \\[5pt]
\mbox{s.t.} & p_i \, = \, p_i^0, & \forall \, i \in \alpha^\infty,\\
&p_i \, \geq  \, p_i^0+\delta_i, & \forall \, i \in \beta^\infty, \\
&p_i \, \leq \, p_i^0-\delta_i, & \forall \, i \in \gamma^\infty.
\end{array}
\end{equation}
Since the objective function in (\ref{eq:restricted_opt}) is strongly convex (due to the positive definiteness of matrix $S$), and all the constraints are linear, there exists a unique optimal solution to (\ref{eq:restricted_opt}). By introducing a multiplier $\mu_i$ for each constraint in (\ref{eq:restricted_opt}), define the Lagrangian function:
$$\mathcal L(\boldsymbol p, \boldsymbol \mu)\, \triangleq \, \mathcal Q(\boldsymbol p) + \sum_{i\in \alpha^\infty} \mu_i (p_i - p_i^0) - \sum_{i\in \beta^\infty} \mu_i (p_i - p_i^0 - \delta_i) + \sum_{i\in \gamma^\infty} \mu_i (p_i - p_i^0+\delta_i).$$
It is well known that a vector $\overline{ \boldsymbol p}$ is the optimal solution of (\ref{eq:restricted_opt}) if and only if there exists a vector $\boldsymbol \mu$ so that the pair $(\overline{\boldsymbol p} , \boldsymbol \mu)$ satisfies:
\begin{equation}\label{eq:KKT_restricted_opt}
\begin{array}{rcll}
\nabla_i \mathcal Q(\overline{\boldsymbol p}) + \mu_i & = & 0, & \forall \, i \in \alpha^\infty \cup \gamma^\infty \\
\nabla_i \mathcal Q(\overline{\boldsymbol p}) - \mu_i & = & 0, & \forall \, i \in \beta^\infty \\
\overline p_i - p_i^0 & = & 0, & \forall \, i \in \alpha^\infty\\
0 \, \leq \, \overline p_i - p^0_i - \delta_i & \perp & \mu_i \geq 0, & \forall \, i \in \beta^\infty\\
0 \, \leq \, -\overline p_i + p^0_i - \delta_i & \perp & \mu_i \geq 0, & \forall \, i \in \gamma^\infty
\end{array}
\end{equation}
where ``$\perp$" means two vectors are perpendicular, i.e., $\boldsymbol x \perp \boldsymbol y \Leftrightarrow \boldsymbol x^T \boldsymbol y = 0$. The following lemma holds.
\begin{lemma}\label{lm:stationary_local_optimal}
Let $\boldsymbol p^\infty \in \mathbb R^n$ be a first-order stationary point defined by (\ref{eq:stationary}). Then it is the optimal solution of (\ref{eq:restricted_opt}) with the index sets $\alpha^\infty$, $\beta^\infty$, and $\gamma^\infty$ as defined in (\ref{eq:index_set_sationary}).
\end{lemma}
\proof{}
It suffices to construct a vector $\boldsymbol \mu$ so that $(\boldsymbol p^\infty, \boldsymbol \mu)$ satisfies (\ref{eq:KKT_restricted_opt}). For each $i \in \alpha^\infty \cup \gamma^\infty$, we let $\mu_i = -\nabla_i \mathcal Q(\boldsymbol p^\infty)$. For each $i\in \beta^\infty$, we let $\mu_i = \nabla_i \mathcal Q(\boldsymbol p^\infty)$. It suffices to verify that the last two complementarity conditions in (\ref{eq:KKT_restricted_opt}). In fact, for any $i \in \beta^\infty$, there are two cases as follows:
\begin{itemize}[leftmargin=.1in]
  \item[] Case 1. $p^\infty_i - p_i^0 - \delta_i>0$. In this case, it is clear that $p^\infty_i \neq p^0_i$. By Lemma \ref{lm:projection_omega_construction_necessary}, we have
  $$ p^\infty_i \in \Pi_{\pdi}\left(p^\infty_i-\frac{1}{L} \nabla_i \mathcal Q(\boldsymbol p^\infty)\right), \,\, \forall \, i \notin \alpha^\infty, $$
  which implies that $p^\infty_i = p^\infty_i-\frac{1}{L} \nabla_i \mathcal Q(\boldsymbol p^\infty) \, \Rightarrow \, \nabla_i \mathcal Q(\boldsymbol p^\infty)=0 .$ Hence, $\mu_i = 0$, and
  ${0 \, \leq \, p^\infty_i - p^0_i - \delta_i} \, \perp \, \mu_i \geq 0$ holds.
  \item[] Case 2. $p^\infty_i - p_i^0 - \delta_i=0$. In this case, we have $p^\infty_i \neq p^0_i$ due to assumption (A4). By Lemma \ref{lm:projection_omega_construction_necessary}, we still have
      $$ p^\infty_i \in \Pi_{\pdi}\left(p^\infty_i-\frac{1}{L} \nabla_i \mathcal Q(\boldsymbol p^\infty)\right), \,\, \forall \, i \notin \alpha^\infty. $$
      By equation (\ref{eq:1d_projection_closedform}) we deduce that
      $$ p_i^0 + \frac{1}{2}\delta_i \, \leq \, p^\infty_i-\frac{1}{L} \nabla_i \mathcal Q(\boldsymbol p^\infty) \, \leq \, p_i^0 + \delta_i\, \Rightarrow \, p^\infty_i-\frac{1}{L} \nabla_i \mathcal Q(\boldsymbol p^\infty) \leq p^\infty_i \, \Rightarrow \, \nabla_i \mathcal Q(\boldsymbol p^\infty) \geq 0.$$
      Thus $\mu_i = \nabla_i \mathcal Q(\boldsymbol p^\infty) \geq 0$, and $0 \, \leq \,  p^\infty_i - p^0_i - \delta_i \, \perp \, \mu_i \geq 0$ holds.
\end{itemize}
We can similarly show that the second complementarity condition in (\ref{eq:KKT_restricted_opt}) holds. This concludes the proof.
 \endproof
Lemma \ref{lm:stationary_local_optimal} indicates that a first-order stationary point satisfying (\ref{eq:stationary}) must be the optimal solution of (\ref{eq:restricted_opt_anypartition}) with the partition being $(\alpha^\infty,\beta^\infty,\gamma^\infty)$ and hence a local optimal solution of (\ref{eq:price_opt_const_base_linear_demand}). On the other hand, not every local optimal solution of (\ref{eq:price_opt_const_base_linear_demand}) is a first-order stationary point. As we have discussed in Section \ref{sec:problem}, each partition $(\alpha, \beta,\gamma) \in \Im$ corresponds to a local optimal solution. We have the following example of a local optimal solution that is not a first-order stationary point.
\begin{example}
Consider the case where $n=2$ and $k=1$, $\boldsymbol p^0 =[0,0]^T$, $\boldsymbol \delta =[0.5,0.5]^T$, $\boldsymbol S = \begin{bmatrix}
	2 & -a \\
	-a & 2
	\end{bmatrix} $, $\boldsymbol r =[-6,-1]^T$, and the problem is
\begin{equation} \label{eq:example}
\begin{array}{rll}
\displaystyle{ \min_{\boldsymbol p} }&  \displaystyle{ \frac{1}{2}\boldsymbol p^T
	\boldsymbol S \boldsymbol p } + \boldsymbol r^T \boldsymbol p. \\
\end{array}
\end{equation}
We can show that the local optimal solution with $\alpha=\{1\}$, $\beta =\{2\}$ and $\gamma = \emptyset$ is not a stationary point. In fact, the local optimal solution with partition $(\{1\},\{2\},\emptyset )$ is $[0,0.5]^T$. Then, we can show that $[0,0.5]^T$ is not a stationary point when $ 0 \leq a < 1$. In fact, when $\boldsymbol p  =[0,0.5]^T$, and $L=3$ we can calculate $\boldsymbol q =\boldsymbol p- \frac{1}{L} \nabla \mathcal Q(\boldsymbol p)$ as follows:
\begin{equation} \label{eq:example}
\begin{array}{lll}
q_1 &=  p_1 - \frac{1}{3}(s_{11}p_1 + s_{12}p_2 -6) & = \frac{a}{3}p_2+  2 >2 \\
  q_2 &=  p_2 - \frac{1}{3}( s_{21}p_1+ s_{22}p_2 -1)  &=  \frac{1}{3}p_2+  \frac{1}{3} = 0.5\\
\end{array}
\end{equation}
To calculate $\mathcal H_{\kpd}(\boldsymbol q)$, we have:
\begin{equation} \label{eq:example}
\begin{array}{lll}
\displaystyle{  \Delta_1 }&= \displaystyle{ (p_1^0 - q_1)^2 - d(q_1, \mathcal P^1_{\delta_1}) } & =  q_1^2\\
\displaystyle{  \Delta_2 }&= \displaystyle{ (p_2^0 - q_2)^2 - d(q_2, \mathcal P^2_{\delta_2}) } & = 0.25\\
\end{array}
\end{equation}
Clearly, $\Delta_1 >\Delta_2$ and hence $[0,0.5]^T \notin \mathcal H_{\kpd}(\boldsymbol q) = \{[q_1, 0]^T\} $. Therefore, the local optimal solution with partition $(\{1\},\{2\},\emptyset)$ is not a stationary point. 
\end{example}
While not every local optimal solution of (\ref{eq:price_opt_const_base_linear_demand}) is a first-order stationary point, the following proposition guarantees that a global optimal solution must be a first-order stationary point.
\begin{proposition}\label{prop:global_opt_stationary}
Let $\boldsymbol p^*$ be a global optimal solution of (\ref{eq:price_opt_const_base_linear_demand}), then $\boldsymbol p^*$ is a first order stationary point.
\end{proposition}
\proof{}
By (\ref{eq:decreasing}) in Lemma \ref{lm:H_property}, any $\boldsymbol \eta \in \mathcal H_{\kpd}\left( \boldsymbol p^* - \frac{1}{L} \nabla \mathcal Q(\boldsymbol p^*) \right)$ satisfies
$$ \mathcal Q(\boldsymbol p^*) - \mathcal Q(\boldsymbol \eta) \, \geq \, \frac{L-\lambda_1}{2} \|p^* - \eta\|_2^2.$$
Since $\boldsymbol p^*$ is a global optimal solution of (\ref{eq:price_opt_const_base_linear_demand}), we must have $ \mathcal Q(\boldsymbol p^*) - \mathcal Q(\boldsymbol \eta) \leq 0$, therefore $\|\boldsymbol p^* - \boldsymbol \eta\|_2^2 = 0$. Hence $\boldsymbol p^* = \boldsymbol \eta \in \mathcal H_{\kpd}\left( \boldsymbol p^* - \frac{1}{L} \nabla \mathcal Q(\boldsymbol p^*)\right)$. This concludes the proof.
 \endproof
In summary, we have shown the following relationship regarding (\ref{eq:price_opt_const_base_linear_demand}):
$$ \{\mbox{global optimal solutions}\} \subseteq \{\mbox{first-order stationary points}\} \subseteq  \{\mbox{local optimal solutions}\}. $$

We next study the quality of a first-order stationary point in terms of the objective value. Let
$\Delta^{\infty}_i \, \triangleq \, (p_i^0 - q_i^\infty)^2- d(q_i^\infty, \pdi),$
with $\boldsymbol q^\infty \, \triangleq \, \boldsymbol p^\infty - \displaystyle{\frac{1}{L} } \nabla \mathcal Q(\boldsymbol p^\infty)$. Assume $\Delta^\infty_{(1)} \geq \Delta^\infty_{(2)} \geq \cdots \geq \Delta^\infty_{(n)}$. The next result provides a performance bound of a first-order stationary point.
\begin{proposition}\label{prop:near_Optimality_general}
Let $\boldsymbol p^\infty$ satisfy the first-order stationary condition \eqref{eq:stationary}, let $\sigma^\infty \, \triangleq \, \{i\, | \, p^\infty_i \neq p_i^0\}$ and $\kappa = |\sigma^\infty|$, then the following two statements hold.
\begin{itemize}[leftmargin=.5in]
\item[(i).] $$\|\nabla \mathcal Q(\boldsymbol p^\infty)\|_2^2 \, \leq \, L^2\sum_{i=\kappa+1}^n\Delta^\infty_{(i)} + \frac{1}{4}  L^2 \sum_{i=1}^n \delta_i^2 .$$
\item[(ii).]   $$\mathcal Q(\boldsymbol p^\infty) - \mathcal Q^* \, \leq \, \frac{L^2}{2\lambda_n}\sum_{i=\kappa+1}^n\Delta^\infty_{(i)} +  \frac{L^2}{8\lambda_n} \sum_{i=1}^n \delta_i^2.  $$
\end{itemize}
\end{proposition}
\proof{}
By definition, for each $i \notin \sigma^\infty$, we have $p_i^{\infty} = p_i^0$. Hence, for each $i \notin\sigma^\infty$
\begin{eqnarray}
&&\Delta_i^{\infty} \, = \, \left(p_i^0 - p_i^\infty -\frac{1}{L} \nabla_i \mathcal Q(\boldsymbol p^\infty)\right)^2 - d(q_i, \pdi) \, = \, \frac{1}{L^2} \left( \nabla_i \mathcal Q(\boldsymbol p^\infty)\right)^2 - d(q_i, \pdi) \nonumber \\
  &\Rightarrow& \left( \nabla_i \mathcal Q(\boldsymbol p^\infty)\right)^2 \, = \, L^2 \left( \Delta_i^\infty + d(q_i, \pdi) \right) \, \leq \, L^2 \left( \Delta_i^\infty + \frac{\delta^2_i}{4} \right),\nonumber
\end{eqnarray}
where the inequality is due to inequality (\ref{eq:distance_upperbound}). For any $i \in \sigma^\infty$, by Lemma \ref{lm:projection_omega_construction_necessary} it is clear that $$p_i^\infty \, \in \, \Pi_{\pdi}\left( p^\infty_i - \frac{1}{L}\nabla_i \mathcal Q(\boldsymbol p^\infty) \right).$$
Notice that for $p_i^\infty \in \pdi$, there are five cases to consider:
\begin{itemize}
  \item[]Case 1. $p_i^\infty > p_i^0 + \delta_i$. In this case, we have
  $$p_i^\infty \, = \, p_i^\infty - \frac{1}{L} \nabla_i \mathcal Q(\boldsymbol p^\infty) \, \Rightarrow \, \nabla_i \mathcal Q(\boldsymbol p^\infty) \,= \,0.$$
  \item[]Case 2. $p_i^\infty = p_i^0 + \delta_i$. In this case, we have
  $$p_i^0 + \frac{1}{2} \delta_i \,\leq \, p_i^\infty-\frac{1}{L} \nabla_i \mathcal Q(\boldsymbol p^\infty) \, \leq \, p_i^0 + \delta_i \, \Rightarrow \, p_i^0 + \frac{1}{2} \delta_i \,\leq \, p_i^0 + \delta_i-\frac{1}{L} \nabla_i \mathcal Q(\boldsymbol p^\infty) \, \leq \, p_i^0 + \delta_i. $$
  Therefore, we have:
  $$ 0 \, \leq \, \nabla_i \mathcal Q(\boldsymbol p^\infty) \, \leq \,\frac{1}{2} L \delta_i.$$
  and hence $(\nabla_i \mathcal Q(\boldsymbol p^\infty))^2 \leq \frac{1}{4} L^2 \delta_i^2$.
  \item[]Case 3. $p_i^\infty = p_i^0$. In this case, $p_i^0 \in \Pi_{\pdi}(q_i)$, which, by equation (\ref{eq:1d_projection_closedform}), implies that
      $$p_i^0 - \frac{1}{2} \delta_i \, \leq \, p_i^\infty - \frac{1}{L} \nabla_i \mathcal Q(\boldsymbol p^\infty) \, \leq \, p_i^0 + \frac{1}{2} \delta_i \, \Rightarrow \, (\nabla_i \mathcal Q(\boldsymbol p^\infty))^2 \, \leq \, \frac{1}{4} L^2 \delta_i^2.$$
  \item[]Case 4. $p_i^\infty < p_i^0 - \delta_i$. Similar to Case 1, we have $\nabla_i \mathcal Q(\boldsymbol p^\infty) = 0$.
  \item[]Case 5. $p_i^\infty = p_i^0 - \delta_i$. Similar to Case 2, we have
  $$ -\frac{1}{2} L \delta_i \, \leq \, \nabla_i \mathcal Q(\boldsymbol p^\infty) \, \leq \, 0 \, \Rightarrow \, (\nabla_i \mathcal Q(\boldsymbol p^\infty))^2 \, \leq \, \frac{1}{4} L^2 \delta_i^2.$$
  \end{itemize}
Overall, for all $i \in \sigma^\infty$, $$(\nabla_i \mathcal Q(\boldsymbol p^\infty))^2 \leq \frac{1}{4}L^2\delta_i^2.$$ Hence, we have
\begin{eqnarray}
\| \nabla \mathcal Q(\boldsymbol p^\infty)\|_2^2 \, = \, \sum_{i=1}^n (\nabla_i \mathcal Q(\boldsymbol p^\infty))^2 & \leq &  \sum_{i \in \{(1),\cdots,(\kappa)\}} \frac{1}{4}L^2\delta_i^2 +  \sum_{i \in \{(\kappa+1),\cdots,(n)\}} L^2\left(\Delta_i^\infty + \frac{1}{4}\delta_i^2\right)\nonumber \\
& = & L^2\sum_{i=\kappa+1}^n\Delta^\infty_{(i)} + \frac{1}{4}  L^2 \sum_{i=1}^n \delta_i^2, \nonumber
\end{eqnarray}
and (i) holds readily. To show (ii), we notice that $\mathcal Q(\boldsymbol p)  =  \frac{1}{2} \boldsymbol p^T S \boldsymbol p - \boldsymbol f^T \boldsymbol p$, with $\boldsymbol f \triangleq \boldsymbol a + \boldsymbol D^T \boldsymbol c$. It is easy to see that $\nabla \mathcal Q(\boldsymbol p) = S \boldsymbol p - \boldsymbol f$. Recall that $\widehat{\boldsymbol p}$ is the unconstrained minimizer of $\mathcal Q(\boldsymbol p)$, and hence $\nabla \mathcal Q(\boldsymbol p) = \boldsymbol 0_n$. Therefore, we have
$$\|S (\boldsymbol p^\infty - \widehat{\boldsymbol p})\|_2 \, = \, \|\nabla \mathcal Q(\boldsymbol p^\infty) - \nabla \mathcal Q(\widehat{\boldsymbol p})\|_2 \, = \, \|\nabla \mathcal Q(\boldsymbol p^\infty)\|_2. $$
By Lemma \ref{lm:linear_demand_properties} item (b), we know that
$$\mathcal Q(\boldsymbol p^\infty) - \widehat{\mathcal Q} \, \leq \, \frac{1}{2\lambda_n}\|S(\boldsymbol p^\infty-\widehat{\boldsymbol p})\|_2^2 \, = \,\frac{1}{2\lambda_n}\|\nabla \mathcal Q(\boldsymbol p^\infty)\|_2^2 \, \leq \, \frac{L^2}{2\lambda_n}\sum_{i=\kappa+1}^n\Delta^\infty_{(i)} +   \frac{L^2}{8\lambda_n} \sum_{i=1}^n \delta_i^2.  $$
 \endproof
The next result is a corollary to Proposition \ref{prop:near_Optimality_general} when $\boldsymbol p^\infty$ satisfies the zero-norm constraint strictly.
\begin{corollary}\label{col:near_Optimality_strict_zeronorm}
Let $\boldsymbol p^\infty$ satisfy the first-order stationary condition \eqref{eq:stationary} and $\sigma^\infty$ and $\kappa$ be defined in Proposition \ref{prop:near_Optimality_general}. If $|\sigma^\infty| < k$, then the following two statements hold.
\begin{itemize}
\item[(i).] $$\|\nabla \mathcal Q(\boldsymbol p^\infty)\|_2^2 \, \leq \, \frac{1}{4}L^2 \sum_{i=1}^n \delta_i^2.$$
\item[(ii).]   $$\mathcal Q(\boldsymbol p^\infty) - \mathcal Q^* \, \leq \, \frac{ L^2}{8 \lambda_n} \sum_{i=1}^{n}\delta_i^2.  $$
\end{itemize}
\end{corollary}
\proof{}
For each $i \in \alpha^\infty$, by Lemma \ref{lm:projection_omega_construction_necessary} and the assumption that $|\sigma^\infty| < k$, we deduce that $\Delta^\infty_{(i)} = 0$ for all $i=\kappa+1,\cdots,n$. The corollary therefore follows readily.
 \endproof
\begin{remark}
In case when $\delta_i = 0$ for all $i=1,\cdots,n$, Corollary \ref{col:near_Optimality_strict_zeronorm} indicates that if $\|\boldsymbol p^\infty - \boldsymbol p^0\|_0 < k$ then $\boldsymbol p^\infty$ is a global optimal solution of (\ref{eq:price_opt_const_base_linear_demand}). This is an extension of Proposition 5 of \cite{bertsimas2016best}.
\end{remark}
We next study the convergence of the GPA. We first state a technical lemma regarding the decrease in the objective value from iteration to iteration.
\begin{lemma}\label{lm:obj_decrease}
Let $\{\boldsymbol p^t\}_{t = 1, 2, \cdots}$  be a sequence generated by the GPA, for any $L > \lambda_1$. Then the sequence $\{\mathcal Q(\boldsymbol p^t)\}_{t=1,2,\cdots}$ satisfies
\begin{equation}\label{eq:obj_decrease}
\mathcal Q(\boldsymbol p^t)  - \mathcal Q(\boldsymbol p^{t+1})) \, \geq  \, \frac{L- \lambda_1}{2} ||\boldsymbol p^{t+1}- \boldsymbol p^{t}||_2^2, \,\, \forall \, t=1,\cdots.	
\end{equation}
Moveover, $$\lim_{t \rightarrow \infty} \| \boldsymbol p^{t+1}-\boldsymbol p^{t}\|_2 \, \rightarrow \, 0.$$
\end{lemma}
\proof{}
By equation (\ref{eq:decreasing}) in Lemma \ref{lm:H_property} it is clear that
\begin{equation}
\mathcal Q(\boldsymbol p^t)  - \mathcal Q(\boldsymbol p^{t+1}) \, \geq  \, \frac{L- \lambda_1}{2} ||\boldsymbol p^{t+1}- \boldsymbol p^{t}||_2^2, \,\, \forall \, t=1,\cdots.
\end{equation}
Therefore the sequence $\{\mathcal Q(\boldsymbol p^t)\}_{t=1,2,\cdots}$ is monotonically decreasing and is bounded from below by $\widehat{\mathcal Q}$, and hence is convergent, i.e.,
$$\lim_{t\rightarrow \infty} \mathcal Q(\boldsymbol p^t) \, = \, \mathcal Q^\infty.$$
Hence, $\lim_{t\rightarrow \infty} \|\boldsymbol p^{t+1}-\boldsymbol p^t\| = 0.$
 \endproof
Next we prove a result regarding the stability of certain key index sets when the number of iterations is large enough. Let $\{\boldsymbol p^t\}_{t=1,2,\cdots}$ be the sequence generated by the GPA. For each $t$, we define
\begin{equation}\label{eq:index_sets}
\alpha^t \triangleq \{i \,| \,p^t_i = p^0_i\},\,\,\beta^t \triangleq \{i \, | \, p^t_i \geq p_i^0 + \delta_i \}, \, \mbox{ and }\, \gamma^t \triangleq \{ i \, | \, p^t_i \leq p_i^0 - \delta_i\}.
\end{equation}
As we can see, since assumption (A4) holds, for each $t=1,2,\cdots$, the index sets $\alpha^t$, $\beta^t$, and $\gamma^t$ are mutually exclusive and satisfy $\alpha^t \cup \gamma^t \cup \gamma^t = \{1,2,\cdots,n\}$.
\begin{lemma}\label{lm:stability_index_sets}
Let $\{\boldsymbol p^t\}_{t=1,2,\cdots}$ be a sequence generated by the GPA. Let the index sets $\alpha^t$, $\beta^t$, and $\gamma^t$ be defined in (\ref{eq:index_sets}). There exists an integer $T$, and mutually exclusive sets $\alpha^T$, $\beta^T$, and $\gamma^T$ such that
$$ \alpha^t = \alpha^T, \,\, \beta^t=\beta^T, \, \mbox{ and } \, \gamma^t = \gamma^T $$ for all $t \geq T$.
\end{lemma}
\proof{}
Let $\delta_{\min} = \min\{\delta_1,\cdots,\delta_n\}$, the smallest among $\delta_i$'s. By Assumption (A4), $\delta_{\min} > 0$. By Lemma \ref{lm:obj_decrease}, there exists an integer $T$ such that $$ \|\boldsymbol p^t - \boldsymbol p^{t+1}\|_2^2 \leq  \frac{\delta_{\min}^2}{2}.$$
For any $t\geq T$, if $\alpha^t \neq \alpha^{t+1}$, then there are two possibilities. (i) There exists an index $i$ such that $i \in \alpha^t$ and $i\notin \alpha^{t+1}$. (ii) There exists an index $i$ such that $i \in \alpha^{t+1}$ and $i\notin \alpha^{t}$. In either case, $|p^t_i - p_i^{t+1}| \geq \delta_i \geq \delta_{\min}$, and hence
$\|\boldsymbol p^t - \boldsymbol p^{t+1}\|_2^2 \geq \delta_{\min}^2$, which is a contradiction. Therefore, for all $t\geq T$, $\alpha^t = \alpha^{t+1} = \cdots $. Similarly, we have $\beta^t = \beta^{t+1} = \cdots $ and $\gamma^t = \gamma^{t+1} = \cdots $ for all $t\geq T$.
 \endproof

\begin{theorem} \label{th:convergence}
Let $\{\boldsymbol p^t\}_{t=1,2,\cdots}$ be a sequence generated by the GPA. It follows that the sequence $\{\boldsymbol p^t\}_{t=1,2,\cdots}$ converges to a first-order stationary point, i.e.,
  $$\lim_{t \rightarrow \infty} \boldsymbol p^t \, = \, \boldsymbol p^\infty,$$
and $\boldsymbol p^\infty$ is a first-order stationary point satisfying (\ref{eq:stationary}).
\end{theorem}
\proof{}
By Lemma \ref{lm:stability_index_sets}, there exists a $T \geq 1$ such that $$ \alpha^t = \alpha^T, \,\, \beta^t=\beta^T, \, \mbox{ and } \, \gamma^t = \gamma^T $$ for all $t \geq T$. Therefore, the GPA is equivalent to applying a gradient projection algorithm on the following convex quadratic program
\begin{equation}\label{eq:stable_QP}
\begin{array}{rll}
\overline{\boldsymbol p}_{\alpha,\beta,\gamma} \, \triangleq \, \displaystyle{ \argmin_{\boldsymbol p} } & \mathcal Q(\boldsymbol p) \\[5pt]
\mbox{s.t.} & p_i \, = \, p_i^0, & \forall \, i \in \alpha^T, \\
&p_i \, \geq  \, p_i^0+\delta_i, & \forall \, i \in \beta^T, \\
&p_i \, \leq \, p_i^0-\delta_i, & \forall \, i \in \gamma^T.
\end{array}
\end{equation}
By the coercivity of $\mathcal Q(\boldsymbol p)$ (part (c) of Lemma \ref{lm:linear_demand_properties}), and Lemma \ref{lm:obj_decrease}, it is clear that $\{\boldsymbol p^i\}_{t=1,2,\cdots}$ is bounded. Hence, the convergence follows \citep{nesterov2013introductory}. Next, we show that the GPA converges to a first-order stationary point satisfying (\ref{eq:stationary}). Because of the continuity of the $\nabla \mathcal Q(\boldsymbol p)$ of a quadratic function,
$$\lim_{t \rightarrow \infty} \left[\boldsymbol p^t- \frac{1}{L} \nabla \mathcal Q(\boldsymbol p^t) \right] \, = \, \boldsymbol p^\infty- \frac{1}{L} \nabla \mathcal Q(\boldsymbol p^\infty).$$
Then, by taking the limit on both sides of (\ref{eq:iter}), because of the closedness of $\varphi(\boldsymbol q)$ by Lemma \ref{lm:closed_map},  we know that $\boldsymbol p^\infty \in \mathcal H_{\kpd}\left(\boldsymbol p^\infty- \frac{1}{L} \nabla \mathcal Q(\boldsymbol p^\infty)\right) $. This completes the proof.
 \endproof
\begin{remark}
From Lemma \ref{lm:stability_index_sets} and Theorem \ref{th:convergence} we can see that once the GPA progresses to the point where $\|\boldsymbol p^{t+1} - \boldsymbol p^t\|_2^2 \leq \frac{\delta^2_{\min}}{2}$, it is equivalent to applying a gradient projection algorithm to the convex quadratic program (\ref{eq:stable_QP}). Since the gradient projection algorithm is a first-order algorithm which converges at most linearly, we can speed up the convergence from that point by switching to a more efficient convex quadratic programming algorithm, such as an interior point algorithm.
\end{remark}

\section{Numerical Experiments} \label{sec:numerical_results}
We conduct numerical experiments on both randomly generated data as well as real-world data. On randomly generated data, we conduct two different groups of numerical experiments. In the first group, we consider the base formulation (\ref{eq:price_opt_const_base_linear_demand}). We include upper/lower bounds on the individual prices in the second group. For both groups and the real-world data, we demonstrate the performance of the GPA by comparing its results to those from Gurobi. 

\subsection{Randomly Generated Data Sets}
Given that the number of products in a typical grocery store in North America is about 50,000 and about 20,000 in Europe, we generate five sets of data corresponding to 10,000, 25,000, 50,000, 75,000 and 100,000 products as follows:
\begin{itemize}[leftmargin=.5in]
\item[1.] Generate matrix $\boldsymbol D$. For each row, we first generate a uniformly distributed random number between $[1,10]$ as the diagonal element, and then generate up to five negative  off-diagonal elements of magnitudes less than 0.2 times that of the diagonal element in that row.
\item[2.] Generate vectors $\boldsymbol a$, $\boldsymbol c$, $\boldsymbol \delta$ and $\boldsymbol p^0$.
    For each product $i$, $p^0_i$ is set to be a random number between $[1,10]$. The vector $-\boldsymbol a - \boldsymbol D^T \boldsymbol c$ is generated as a random vector with all elements in $[1,10]$. $\boldsymbol \delta$ is set as $\boldsymbol 1_n$ for \$1.00 and $\boldsymbol \delta = 0.5 \cdot \boldsymbol 1_n$ for \$0.50.
\item[3.] Set $k$. We allow up to 10\% of the products to change prices.
\end{itemize}
Counting the number of data sets and the two possibilities of $\boldsymbol \delta$ for each case, we have a total of 10 instances of the price optimization problem. For each instance, we use five distinct initial solutions to start the GPA: (1) $\boldsymbol p^0$, (2),(3),(4) are three random vectors, and (5) the solution obtained by using $\boldsymbol p^0$ with a step length larger than $\frac{1}{L}$ for some $L$ greater than the largest eigenvalue of $\boldsymbol D$. The motivation of using initial solution (5) is that by allowing the algorithm to take a longer step, we hope that it can move away from the local optimal solution. In practice, we observe that the best solution often comes from (5). Together with the five initial solutions, we have a total of 50 optimization problems that are solved sequentially.

\subsection{Benchmark Formulation}
We use the results from Gurobi as a benchmark to demonstrate the performance of the GPA. To solve the price optimization (\ref{eq:price_opt_const_base_linear_demand}) using Gurobi, for each product $i$ we introduce a triplet of binary variables, $z_i^P$, $z_i^R$, $z_i^L$ indicating if the price of product $i$ stays at $p^0_i$, is at least $p_i^0+\delta_i$, or is up to $p_i^0 - \delta_i$, respectively. We use the following mixed integer programming formulation:
\begin{equation}\label{eq:MIP_formulation}
\begin{array}{rlll}
\max_{\boldsymbol p, \boldsymbol z^P, \boldsymbol z^R, \boldsymbol z^L} \quad & \displaystyle{- \frac{1}{2}\boldsymbol p^T S \boldsymbol p + (\boldsymbol a + \boldsymbol D^T \boldsymbol c)^T \boldsymbol p } \\[5pt]
\textrm{s.t.} \quad & \displaystyle { p_i \geq p_i^0 z_i^P -M z_i^L + (p_i^0+\delta_i) z_i^R } ,& \forall \, i = 1,\cdots,n, & -- \text{(lower bound)}\\[5pt]
&\displaystyle{ p_i \leq p_i^0 z_i^P + (p_i^0-\delta_i) z_i^L + M z_i^R } ,& \forall \, i = 1,\cdots,n, & -- \text{(upper bound)}\\[5pt]
&\displaystyle{ z_i^P + z_i^L +  z_i^R =1 } ,& \forall \, i = 1,\cdots,n, \\[5pt]
& \displaystyle{ \sum_{i =1}^n (z_i^L + z_i^R) } \leq k, \\[5pt]
& \displaystyle{ z_i^P, z_i^L, z_i^R \in \{0,1\}}, & \forall \, i = 1,\cdots,n,
\end{array}
\end{equation}
where $M$ is a large positive integer chosen to be relatively small to help with the performance of Gurobi. When we have lower/upper bounds $l_i \leq p^0_i-\delta_i$ and $u_i \geq p^0_i +\delta_i$ on the individual prices $p_i$, we use the following formulation, which is similar to the above formulation with $-M$ and $M$ being replaced by the lower bounds and upper bounds, respectively.
\begin{equation}\label{eq:MIP_formulation_bounds}
\begin{array}{rlll}
\max_{\boldsymbol p, \boldsymbol z^P, \boldsymbol z^R, \boldsymbol z^L} \quad & \displaystyle{ -\frac{1}{2}\boldsymbol p^T S \boldsymbol p + (\boldsymbol a + \boldsymbol D^T \boldsymbol c)^T \boldsymbol p } \\[5pt]
\textrm{s.t.} \quad
&\displaystyle{ p_i \geq p_i^0 z_i^P + l_i z_i^L + (p_i^0+\delta_i) z_i^R }, & \forall \, i =1,\cdots,n, & -- \text{(lower bound)}\\[5pt]
&\displaystyle{ p_i \leq p_i^0 z_i^P + (p_i^0-\delta_i) z_i^L + u_i z_i^R }, & \forall \, i =1,\cdots,n, & -- \text{(upper bound)}\\[5pt]
&\displaystyle{ z_i^P + z_i^L +  z_i^R =1} ,& \forall \, i =1,\cdots,n, \\[5pt]
& \displaystyle{ \sum_{i=1}^n (z_i^L + z_i^R) \leq k }, \\[5pt]
&\displaystyle{ z_i^P, z_i^L, z_i^R \in \{0,1\} },& \forall \, i=1,\cdots,n.
\end{array}
\end{equation}
To compare the solution quality, we used the adjusted objective gap as the performance measure:
$$ \text{Adjusted Objective Gap} \, \triangleq \, \frac{100(Z(\boldsymbol p^B) - Z(\boldsymbol p^G))}{|Z(\boldsymbol p^0)|},$$ 
where $\boldsymbol p^0$ is the baseline price, $\boldsymbol p^G$ is the Gurobi solution and $\boldsymbol p^B$ is the GPA solution.  When the adjusted objective gap is positive, the GPA achieved a better solution than Gurobi, whereas when it is negative, Gurobi obtained a better solution than the GPA. We use Gurobi 7.0.2 configured as follows: the maximum run time is set to 1 hour for experiments labeled Gurobi-1 and 4 hours for Gurobi-4, the acceptable optimality gap is set to 0, and all other parameters are at their default values. The number of computing threads is set to 1 to allow for a fair comparison with the GPA, which uses only 1 thread. Both Gurobi and GPA implementations were executed on Amazon Elastic Computing Cloud running Amazon Linux AMI 2018.03 with dual Intel\textsuperscript{\textregistered} Xeon\textsuperscript{\textregistered} 16-core E5-2686 v4 CPU @ 2.3GHz and 244GB memory.

\subsection{Results from the Randomly Generated Data Sets}
\begin{table}
	\caption{Adjusted objective gaps on random data sets.}
	\begin{center}
	    \begin{tabular}{|l|r|rrc|rrc|}
		 	\hline
		 	\multirow{2}{*}{\hspace{0.14in}Bounds}&\multirow{2}{*}{$n$\hspace{0.2in}}&\multicolumn{3}{c|}{Gurobi-1}&\multicolumn{3}{c|}{Gurobi-4}\\
		 	\cline{3-8}
		 	& & Mean & StDev & Range & Mean & StDev & Range \\
			\hline
			\multirow{5}{*}{no bounds}
			& 10,000 & -3.91 & 3.97 & [-6.24, 4.62] & -5.86 & 0.34 & [-6.25, -5.39] \\
			& 25,000 & 29.37& 0.56 & [28.53, 30.10] & 29.37 & 0.56 & [28.53, 30.10] \\
			& 50,000 & 29.58 & 0.70 & [28.26, 30.49] & 29.58 & 0.70 & [28.26, 30.49] \\
			& 75,000 & 29.19 & 0.52 & [28.46, 30.08] & 29.19 & 0.52 & [28.46, 30.08] \\
			& 100,000 & 29.25 & 0.44 & [28.68, 30.14] & 29.25 & 0.44 & [28.68, 30.14] \\
			\hline
			\multirow{5}{*}{\shortstack{$l_i \in [1,5]$ \\ $u_i \in [5,10]$}}
			& 10,000 & -2.19 & 0.18 & [-2.48, -1.99] & -2.19 & 0.18 & [-2.48, -2.00] \\
			& 25,000 & -2.26 & 0.11 & [-2.43, -2.10] & -2.26 & 0.11 & [-2.43, -2.11] \\
			& 50,000 & 22.89 & 13.26 & [-2.33, 29.51] & -2.25 & 0.06 & [-2.35, -2.17] \\
			& 75,000 & 28.86 & 0.24 & [28.39, 29.16] & 4.08 & 13.12 & [-2.36, 29.16] \\
			& 100,000 & 28.86 & 0.25 & [28.48, 29.29] & 28.86 & 0.25 & [28.48, 29.29] \\
			\hline
			\multirow{5}{*}{\shortstack{ $l_i \in [1,5]$ \\ $u_i \in [10,15]$}}
			& 10,000 & -4.00 & 0.11 & [-4.15, -3.89] & -4.02 & 0.11 & [-4.16, -3.91] \\
			& 25,000 & -3.76 & 0.18 & [-3.96, -3.35] & -4.13 & 0.04 & [-4.18, -4.05] \\
			& 50,000 & 26.62 & 0.43 & [26.01, 27.32] & 26.62 & 0.43 & [26.01, 27.32] \\
			& 75,000 & 26.48 & 0.26 & [26.00, 26.85] & 26.48 & 0.26 & [26.00, 26.85] \\
			& 100,000 & 26.51 & 0.28 & [25.91, 26.94] & 26.51 & 0.28 & [25.91, 26.94] \\
			\hline
			\multirow{5}{*}{\shortstack{ $l_i \in [1,5]$ \\ $u_i \in [15,20]$}}
			& 10,000 & -4.33 & 0.19 & [-4.56, -3.94] & -4.39 & 0.14 & [-4.58, -4.20] \\
			& 25,000 & 20.26 & 12.09 & [-2.85, 26.35] & 1.86 & 12.64 & [-4.43, 25.93] \\
			& 50,000 & 26.13 & 0.35 & [25.43, 26.57] & 26.13 & 0.35 & [25.43, 26.57] \\
			& 75,000 & 26.08 & 0.29 & [25.67, 26.56] & 26.08 & 0.29 & [25.67, 26.56] \\
			& 100,000 & 26.05 & 0.22 & [25.69, 26.38] & 26.05 & 0.22 & [25.69, 26.38] \\
			\hline
		\end{tabular}
	\end{center}
	\label{tb:random_data_results}
\end{table}

For a reasonable total run time of 205 minutes solving the 50 optimization problems for the unbounded case, Table~\ref{tb:random_data_results} reports on the combined results for cases when the minimum change is \$0.50 and \$1.00. The results are combined as they are comparable. The GPA significantly outperforms both Gurobi-1 and Gurobi-4 when $n > 10,000$, leading to incremental profits of around 29\% from the baseline. In fact, Gurobi barely improved the objective value from 1 hour to 4 hours, achieving about 35\% MIP gap for $n = 10,000$ but no meaningful MIP gaps for larger cases at the end of 4 hours. For $n = 10,000$, the average adjusted objective gap compared with Gurobi-1 is $-3.91\%$ with a maximum of $-6.24\%$, and the average adjusted objecive gap compared with Gurobi-4 is $-5.86\%$ with a maximum of $-6.25\%$. These results indicate that the GPA is better suited for real-world applications, especially when the number of products is large.

Next, we report on experiments with upper/lower bounds on the prices of individual products. For each product, we generate a random lower bound between $l_i \in[1,5]$ and three increasing random upper bounds within $u_i \in [5,10]$, $u_i \in [10,15]$ and $u_i \in [15,20]$. $p_i^0$ is set to be a random number between $[l_i,u_i]$. The purpose of considering three different upper bounds is to investigate the impact of the width of $[l_i, u_i]$ on the performance of the GPA.

As expected, the bounds significantly reduced the feasible regions of the problem instances, and so both the GPA and Gurobi solved the problems faster. The GPA took less time to converge. For example, with $u_i\in[5,10]$, the GPA took up to 19 minutes on the 50 problems, whereas Gurobi-1 and Gurobi-4 ended up with smaller MIP gaps (6\%) for the smaller cases ($n <= 25,000$ for Gurobi-1 and $n <= 50,000$ for Gurobi-4). As the upper bounds increase, the running time of the GPA increases, i.e., up to 23 minutes total run time for $u_i \in [10,15]$ and up to 30 minutes for $u_i \in [15,20]$. The maximum run times are all for $n = 100,000$. 

Similar to the unbounded case, the GPA performed better than Gurobi on larger instances on all three different upper bounds. For  example when $u_i \in [5,10]$, Table~\ref{tb:random_data_results} shows that the GPA outperforms Gurobi-1 when $n \geq 50,000$, and outperforms Gurobi-4 when $n \geq 75,000$. The critical number of products at which GPA outperforms Gurobi decreases as the upper bounds increase, i.e., the GPA performs better than Gurobi with a larger feasible region.

It is also worth noting that for the instances where Gurobi outperforms the GPA, its advantage as measured by the adjusted objective gap generally decreases with the feasible region. This means that for smaller problem instances with reasonably tight bounds, the GPA can reliably provide good quality solutions in significantly shorter time.

\subsection{Numerical Experiments on Real-World Cases}
\label{sec:real}
We conduct experiments on real-world data sets from a grocery chain with 300 stores in four regions in Europe. Based on a rule where the prices of each product is set uniformly across stores in each region but may differ across regions, we have 59,000 prices to optimize.  We also consider both substitute and complement products, i.e., positive and negative off-diagonal elements in matrix $\boldsymbol D$. This relaxes assumption (A1), but the matrix $\boldsymbol S$ remains positive definite.

First, we consider three cases of fixed minimum price change thresholds across all products: 1 cent, 10 cents and 50 cents. We also consider a cases of varying minimum price changes set to  10\% of the baseline price of each products, i.e., $\delta_i = 10\% \times p^0i$. With the four cases of minimum price change, maximum number of prices that can be changed is limited to 3,000, 5,000 and 8,000. In these experiments, the GPA solutions for each case are obtained in about 10 minutes on sequential runs with 5 initial solutions. The results are summarized in Table \ref{tb:real_world_data_results}. In comparison, Gurobi, initialized with the baseline prices, did not provide better solutions than $\boldsymbol p^0$ for cases with 3,000 and 5,000 maximum number of price changes and did not provide better solutions than the GPA in any of the cases within 4 hours.

\begin{table}
	\caption{Results on real-world cases.}
	\begin{center}
		\begin{tabular}{| c c | c c | c |}
			\hline
			\multirow{3}{*}{\shortstack{Min. Amount of \\ Price Change}} & \multirow{3}{*}{\shortstack{Max. No. of \\ Price Changes}} & \multicolumn{2}{c|}{\textbf{GPA}} & \textbf{Gurobi-4} \\
			\cline{3-5}
			& & Improve from Base & Run Time & Improve from Base \\
		    & & (million \euro{}) &(minutes) & (million \euro{}) \\
			\hline
			\multirow{3}{*}{1 cent} & 3,000 & 615 (15.4\%) & 10.3 & 0 \\
						& 5,000 & 696 (17.5\%) & 7.1 & 0 \\
				   		& 8,000 & 777 (19.5\%) & 10.0 & 546 (13.7\%) \\
			\hline
			\multirow{3}{*}{10 cents} & 3,000 & 611 (15.3\%) & 10.1 & 0 \\
						& 5,000 & 693 (17.4\%) & 7.3  & 0 \\
						& 8,000 & 773 (19.4\%) & 10.3 & 567 (14.2\%) \\
			\hline
			\multirow{3}{*}{50 cents} & 3,000 & 569 (14.3\%) & 10.0 & 0 \\
						& 5,000 & 635 (15.9\%) & 9.0 & 0 \\
						& 8,000 & 708 (17.7\%) & 9.3 & 554 (13.9\%) \\
			\hline
			\multirow{3}{*}{10\% of $p_i^0$} & 3,000 & 611 (15.3\%) & 10.0 & 0 \\
						& 5,000 & 691 (17.3\%) & 6.9 & 0 \\
						& 8,000 & 772 (19.3\%) & 9.7 & 560 (14.0\%) \\
			\hline
		\end{tabular}
	\end{center}
\label{tb:real_world_data_results}
\end{table}

\begin{figure}[h]
	\centering
	\includegraphics[scale=.7]{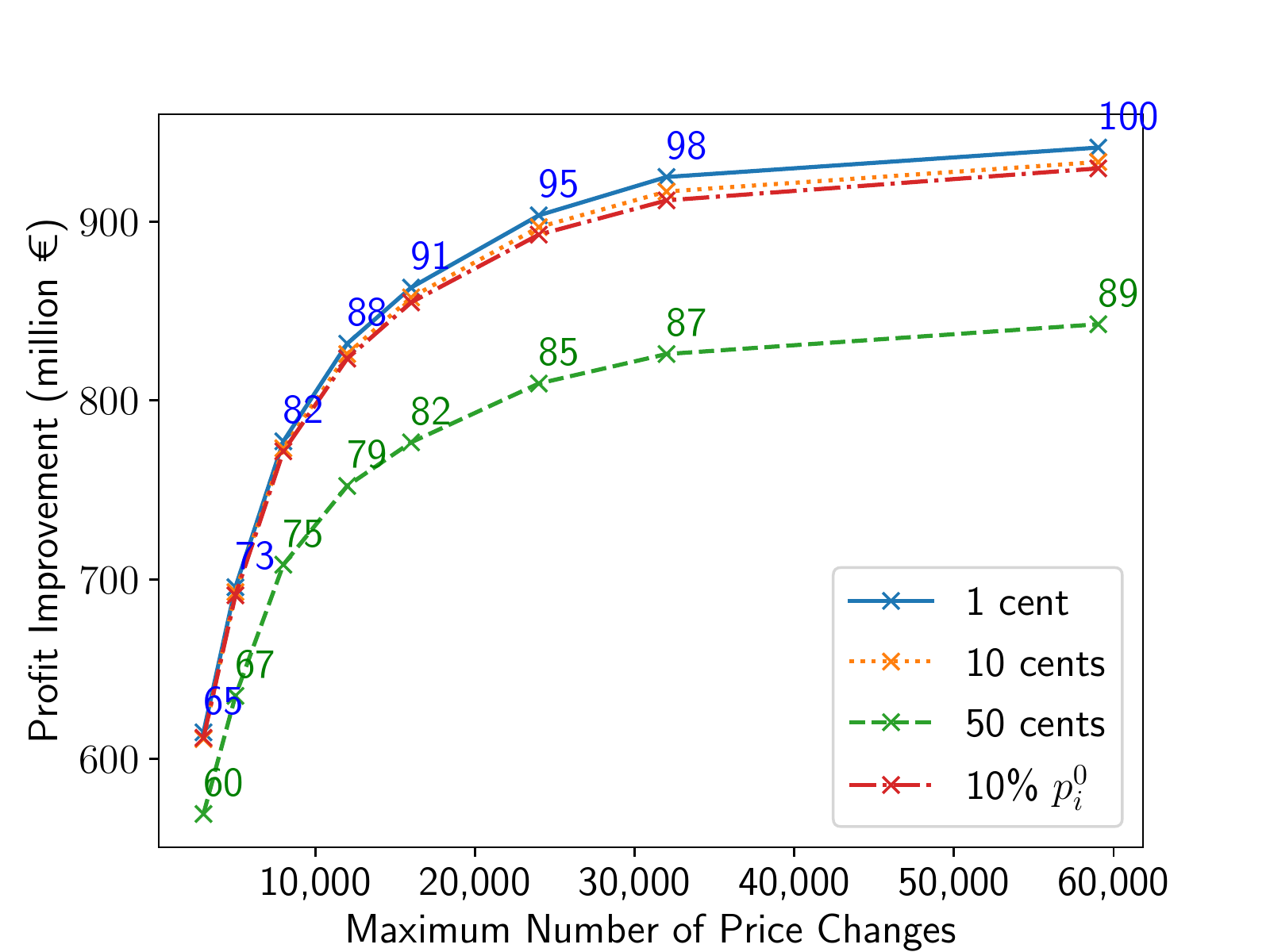}
	\caption{\centering Real-world application of the GPA: Diminishing returns over maximum number of price changes.}
	\label{fig:real}
\end{figure}

We also observed that increasing the maximum number of price changes provide diminishing returns on maximum profit (Figure~\ref{fig:real}). Using as reference the profit improvement obtained by the GPA with minimum price change threshold of 1 cent and $k = 59,000$ (i.e., 100\% of potential profit improvement), we then observe that $k = 3,000$ (5.1\%) can achieve more than 65\% of potential profit improvement, $k = 8,000$ (13.5\%) can achieve more than 80\% potential profit improvement, and  $k = 24,000$ (40.7\%) can achieve more than 95\% potential profit improvement. Notice that since the GPA only finds local optima, it may not be the global maximum improvement. Moreover, a minimum change of 1 cent, 10 cents or 10\% of the baseline price $\boldsymbol p^0$ lead to similar results, while a 50 cent minimum change threshold would reduce the profit by 5-11\%.

\section{Conclusion}\label{sec:conclusion}
In this paper, we study a brick-and-mortar retailer price optimization problem with maximum number of price changes and minimum amount of price change constrains. These two constraints are crucial to provide actionable price recommendations to retailers. Although the feasible region of this optimization problem is not convex, we showed that the Euclidean projection onto it can be obtained very efficiently. Therefore, a gradient projection algorithm is proposed to solve this problem. We demonstrated desired theoretical properties of the algorithm and verified its efficiency using computational experiments.

There are several future research directions stemming from this work. One immediate extension is to include other more standard business rules, such as those that can be written as linear inequalities. With those additional constraints, it is not easy to compute the projection of an arbitrary point to the feasible region, and therefore the GPA in this paper is not directly applicable. However, one can consider decomposition approaches, where the GPA can be applied to solve sub problems. Theoretical properties of such decomposition approaches as well as their empirical performance are worth further investigation.  It is also worth investigating if existing algorithms for solving mathematical programs with cardinality constraints, such as those reported in \cite{burdakov2016mathematical}, can be extended to solve the price optimization problem with the two practical constraints. Another potential extension is to consider other demand functions, such as those resulting from discrete choice models. These demand functions are not necessarily convex, which adds another layer of challenge to this problem.

As we have mentioned before, beyond price optimization problems, the two practical constraints studied are considerably common in various applications. Therefore, we believe that optimization problems with these two practical constraints deserve special attention. Mathematical properties of the feasible set defined by these two constraints as well as sophisticated algorithms taking advantage of those properties can be the focus of future research.

%
%
%




\bibliographystyle{plainnat} 
\bibliography{bib} 



\end{document}